\documentclass[10pt]{article}

\usepackage{amsmath}
\usepackage{amsthm}
\usepackage{amssymb}
\usepackage{amsfonts}
\usepackage{graphicx}
\usepackage{verbatim}
\usepackage{mathrsfs}
\usepackage{latexsym}
\usepackage{graphicx}
\usepackage{color}
\usepackage[hidelinks]{hyperref}
\usepackage{bm}
\usepackage{epstopdf}

\usepackage{hhline}

\usepackage{booktabs}
\usepackage{colortbl}

\usepackage{caption}
\usepackage{subcaption}

\newcommand{\R}{\mathbb{R}}

\newcommand{\Z}{\mathbb{Z}}
\newcommand{\XX}{\mathcal{X}}

\newcommand{\M}{\mathcal{M}}

\newcommand{\g}{\bm{g}}
\renewcommand{\u}{\bm{u}}
\renewcommand{\v}{\bm{v}}

\newcommand{\bu}{\bar{\bm{u}}}
\newcommand{\bv}{\bar{\bm{v}}}
\newcommand{\bw}{\bar{\bm{w}}}
\newcommand{\bs}{\bar{\bm{\sigma}}}

\newcommand{\1}{\bm{1}}
\newcommand{\0}{\bm{0}}

\newcommand{\B}{\mathcal{B}}

\renewcommand{\div}{\text{div}}

\theoremstyle{plain}
\newtheorem{theorem}{Theorem}[section]

\newtheorem{remark}[theorem]{Remark}

\title{Computer-assisted proofs for some nonlinear diffusion problems}

\date{\today}

\begin{document}

\author{ Maxime Breden\thanks{CMAP, \'Ecole Polytechnique, route de Saclay, 91120 Palaiseau, France.}}

\maketitle

\begin{abstract}
In the last three decades, powerful computer-assisted techniques have been developed in order to validate a posteriori numerical solutions of semilinear elliptic problems of the form $\Delta u +f(u,\nabla u) = 0$. By studying a well chosen fixed point problem defined around the numerical solution, these techniques make it possible to prove the existence of a solution in an explicit (and usually small) neighborhood the numerical solution. In this work, we develop a similar approach for a broader class of systems, including nonlinear diffusion terms of the form $\Delta \Phi(u)$. In particular, this enables us to obtain new results about steady states of a cross-diffusion system from population dynamics: the (non-triangular) SKT model. We also revisit the idea of automatic differentiation in the context of computer-assisted proof, and propose an alternative approach based on differential-algebraic equations.
\end{abstract}

\section{Introduction}

\subsection{Context}

Diffusion is a key mechanism in many spatially extended system coming from Physics, Chemistry or Biology. Starting from the prototypical mathematical model used to describe such phenomena that is the heat equation
\begin{align*}
\partial_t u = \Delta u,
\end{align*}
many more general models have been introduced, for instance in order to take into account nonlinear diffusion effects
\begin{align}
\label{eq:NLD}
\partial_t u = \Delta \Phi(u).
\end{align}

Some prime examples which have already been studied extensively in the mathematical literatur are the fast-diffusion equation and the porous medium equation, which correspond to $\Phi(u) = u^m$ with $m<1$ and $m>1$ respectively~\cite{Aro86,Vaz07}. More general nonlinearities are also of interest, for instance $\Phi(u) = \log(u)$ which corresponds to the Ricci flow in two dimensions~\cite{Ham82,TopYin17}. 

Nonlinear diffusion also plays a crucial role for systems, where $u=(u_1,\ldots,u_d)$, especially when the diffusion rate of one component can be affected by the other components. This phenomenon is often referred to as cross-diffusion, and typical mathematical models used in this case are of the form 
\begin{align}
\label{eq:cross-diff}
\partial_t u = \div (A(u) \nabla u).
\end{align}
In some cases, namely when the matrix $A(u)$ is the derivative of some map $\Phi$, these systems can also be written under the form~\eqref{eq:NLD}.

For some of these nonlinear diffusion problems, like the fast-diffusion equation or the porous medium equation, the long-time behavior of the solutions is well understood, see for instance the monograph~\cite{Vaz06} and the references therein. For more general nonlinear diffusion equations, and especially for systems, the situation is more complex, but in some cases entropy method can be used to guarantee at least the existence of global in time solutions, see e.g.~\cite{Jun15}. 

However, in most applications the models also include nonlinear reaction terms, that is,~\eqref{eq:NLD} is replaced by an equation of the form
\begin{align*}
\partial_t u = \Delta \Phi(u) + R(u).
\end{align*}

In this case the dynamics can be way more complex, therefore the study of the long time behavior is usually significantly more complicated. Even in some specific cases where one can prove that the solution converges to a steady state, see e.g. the early work~\cite{LanPhi85}, solving the corresponding stationary problem
\begin{align}
\label{eq:NLeq}
0 = \Delta \Phi(u) + R(u),
\end{align}
is still a demanding task, especially if one is not only interested in existence results but also in quantitative information about the solution(s).

\medskip

In this work, we develop a methodology to get quantitative existence results about equations of the form~\eqref{eq:NLeq} on bounded domains, based on computer-assisted proofs. Before giving more details about the scope and the limitation of our results, let us briefly review some past computer-assisted works on elliptic equations upon which we build.

\subsection{Computer-assisted proofs for elliptic equations}
\label{sec:CAP_Ell}

Computer-assisted proofs for elliptic equations originate from the pioneering works of Nakao~\cite{Nak88} and Plum~\cite{Plu92}, and were then further developed and popularized by them and many others, see e.g.~\cite{Ois95,Yam98,ZglMis01,AriKoc10,DayLesMis07,TakLiuOis13,BerWil17,Wan17,Gom19}.

Most of these works share the same general approach, which consists in validating a posteriori an approximate solution $\bar{u}$ (for an alternative approach, having a topological rather than functional analytic flavor, see~\cite{ZglMis01}). 
In order to conduct this \emph{a posteriori validation}, one considers the elliptic PDE as an $F(u)=0$ problem on some well chosen function space, and then turns it into a fixed point problem $T(u)=u$, where the fixed point operator $T$ should be a contraction\footnote{Historically, some of the early works used a slightly different strategy based on Schauder's fixed point theorem. However, it seems most of the techniques gradually evolved and aligned towards the usage of a contraction mapping, if only because that is often the easiest way to enforce the inclusion assumption needed for Schauder's fixed point theorem.} in a neighborhood of $\bar{u}$. In practice, the two essentials steps are the definition of a suitable fixed point operator $T$, and the derivation of explicit estimates allowing to apply Banach's fixed point theorem to $T$ in a neighborhood of $\bar{u}$. So as to consistently get a contraction, $T$ is usually defined as something close to
\begin{align*}
T(u) = u - DF(\bar u)^{-1} F(u),
\end{align*}
where $DF$ denotes the Frechet derivative of $F$. Provided $\bar u$ is a sufficiently good approximate solution, meaning that $\Vert F(\bar u)\Vert$ is small enough, the key requirement for applying Banach's fixed point theorem to $T$ is to get an explicit control on $\Vert DF(\bar u)^{-1} \Vert$. At this stage, there are three main options.
\begin{enumerate}
\item If the problem is defined on Hilbert spaces, eigenvalues bounds can be used to directly estimate the norm of $DF(\bar u)^{-1}$~\cite[Part II]{NakPluWat19}. While this approach is typically more demanding in practice than the other two presented below, it is also more versatile, and can for instance handle unbounded domains.
\item On bounded domains, an alternative approach consists in studying $DF(\bar u)^{-1}$ by introducing a well chosen \emph{approximate inverse} $A$ of $DF(\bar u)$, which is obtained using a finite dimensional projection. In that case, it is more or less equivalent to directly work with $T$ defined by
\begin{align*}
T(u) = u - A F(u).
\end{align*}
Since $A$ is then defined rather explicitly, computing $\Vert A\Vert$ becomes straightforward, but the difficulty is transferred to proving that $A$ is a good enough approximate inverse, meaning that $\Vert I -ADF(\bar u)\Vert <1$, in order for $T$ to still be a contraction.
\item Still on bounded domains, a third option was recently introduced~\cite{SekNakOis20}. While it aims at directly studying $DF(\bar u)^{-1}$ itself like the first option, it does so by also introducing a finite dimensional projection to decompose $DF(\bar u)$ and then using the Schur complement. Therefore, in practice this approach ends up being close to the second one, as it requires very similar estimates, but those are combined together slightly differently in the end.
\end{enumerate}
Of course, in order to prove that $T$ is a contraction one also has to control the nonlinear terms, but since such control is only required on a small neighborhood of $\bar u$, this is usually not the main difficulty. The whole procedure boils down to a kind of Newton-Kantorovich theorem~\cite{Ort68}, an example of which will be presented in details in Section~\ref{sec:fixed-point}.

As we already mentioned, the general strategy and its variations that we just described were successfully applied to several elliptic problems. However, in all these studies, the leading differential operator is always simply $\Delta u$, or $\Delta^2 u$. The only exceptions seem to be~\cite{BreDevLes15} and~\cite{BreEng21}, where some very specific non-constant coefficients in front of the Laplacian are considered.

\subsection{Generalization to some nonlinear diffusion problems}

In this paper, we generalize the framework presented above, and more specifically option 2, in order to treat problems of the form
\begin{align}
\label{eq:NLeqBC}
\left\{
\begin{aligned}
&\Delta \Phi(u) + R(u) = 0 \qquad & \text{on }\Omega \\
&\frac{\partial u}{\partial n } = 0 \qquad & \text{on } \partial \Omega
\end{aligned}
\right.
\end{align}
on a bounded domains. The key step is to introduce an approximate inverse $A$ with a different structure than those which are typically used for linear diffusion. 

Before going further, let us discuss the assumptions we are going to make in this work, and specify how relevant they are.
\begin{itemize}
\item We restrict our attention to rectangular domains, i.e. of the form $\prod [a_i,b_i]$. This allows us to use a discretization based on Fourier series, which makes the derivation of some of the estimates easier. However, we do not consider this assumption to be essential: the same ideas could be adapted to a discretization based on finite elements, in order to handle more general domains.
\item Similarly, we consider Neumann boundary conditions because they are natural for the main application we have in mind (namely the SKT system which describe the density of several species in a bounded environment), but the analysis would be similar with Dirichlet boundary conditions.
\item To finish with the domain, we rely crucially on the fact that the domain is bounded, and a significantly different approach would be needed for unbounded domains. While we did not explore this possibility thoughtfully, we believe that option 1 described in Section~\ref{sec:CAP_Ell} could be generalized to handle nonlinear diffusion terms as in~\eqref{eq:NLeqBC}.
\item We chose not to include a first order term $\nabla u$ in~\eqref{eq:NLeqBC} in order to keep the presentation as simple as possible, but adding such a term presents no essential difficulty.
\item Still in order to limit technicalities and to focus on the main novelty of this work, namely the treatment of the nonlinear diffusion, we only consider examples where the reaction terms $R$ are at most quadratic. This can also be easily generalized, and we make some remarks in that regard later on.
\item We also emphasize already that our method does not require specific assumptions on the nonlinearity $\Phi$, except smoothness. We do have some non-degeneracy condition, but it will appear as an a posteriori information, rather than as an a priori requirement. That is, if we manage to validate a solution $u$ to~\eqref{eq:NLeqBC}, it will mean that $D\Phi(u)$ is nonsingular, but we do not need to know this in advance for the method to work. For instance, it could be that $D\Phi(x)$ is singular for some values of $x$, but that the solution $u$ ends up never taking these values. In some sense, this means the problem~\eqref{eq:NLeqBC} could have been rewritten as
\begin{align*}
\Delta v + R\left(\Phi^{-1} (v) \right) = 0,
\end{align*} 
but our analysis never requires us to deal with $\Phi^{-1}$. Finally, let us point out that, as soon as we know that $D\Phi(u)$ is nonsingular, we can rewrite the boundary conditions in~\eqref{eq:NLeqBC} in a form which directly involves fluxes, namely $\frac{\partial \Phi(u)}{\partial n } = 0$.
\end{itemize}

\begin{remark}
The approach developed in this work is also suitable for uniformly elliptic problems of the form $\Delta (a u) + R(u) = 0$, where $a$ is a non constant coefficient (or matrix in the case of systems), which already fall outside of the usual framework for computed-assisted proofs for equations of the form  $\Delta u + R(u) = 0$.
\end{remark}

The remainder of this paper is organized as follows. In Section~\ref{sec:PM} we start with a simple case, namely a scalar equation with $\Phi(u)=u^2$, in order to expose the main ideas without to many technicalities. In Section~\ref{sec:SKT} we then focus on the SKT system, and present some new results about its steady states. Finally, we discuss in Section~\ref{sec:NP} how to handle more complicated, non-polynomial, diffusion terms, and propose an alternative to the usual automatic-differentiation technique. 

All the Matlab codes used for the computer-assisted parts of the proofs are available at~\cite{Bre21}.

\section{A simple case first}
\label{sec:PM}

In this section, we consider a scalar problem on a one-dimensional domain of the form

\begin{align}
\label{eq:PM}
\left\{
\begin{aligned}
&\Delta \Phi(u) + R(u) = 0 \qquad & \text{on }(0,1), \\
&\frac{\partial u}{\partial n } = 0 \qquad & \text{on } \{0,1\},
\end{aligned}
\right.
\end{align}
where $\Phi(u)=u^2$ and $R(u) = \alpha u -\beta u^2 + g$.

\subsection{Notations and sequence spaces}
\label{sec:not}

The material presented in this subsection is standard, and mostly included for the sake of fixing notations.

We look for solutions as Fourier series, i.e.
\begin{equation*}
u(x)= u_0 + 2\sum_{n=1}^\infty u_n \cos(n\pi x) = \sum_{n\in\Z} u_{\vert n\vert} \cos(n\pi x),
\end{equation*}
and denote by $\u = \left(u_n\right)_{n\geq 0}$ the sequence of Fourier coefficients associated to a function $u$. In the sequel, we always use this convention that a bold symbol denotes the sequence of Fourier coefficients associated to the function of the same name. For instance $\1$ is the sequence $(1,0,\ldots,0,\ldots)$ representing the constant function equal to 1.

For any $\nu\geq 1$, we consider
\begin{align*}
\left\Vert \u \right\Vert_{\nu} &= \vert u_0\vert + 2\sum_{n=1}^\infty \vert u_n\vert \nu^n \\
&= \sum_{n\geq 0} \vert u_n\vert \xi_n(\nu),
\end{align*}
where
\begin{align*}
\xi_n(\nu) = 
\left\{
\begin{aligned}
&1 \qquad & n=0,\\
&2\nu^n \qquad & n\geq 1,
\end{aligned}
\right.
\end{align*}
and the associated Banach space 
\begin{align*}
\ell^1_\nu = \left\{\u = \left(u_n\right)_{n\geq 0},\  \left\Vert \u \right\Vert_{\nu}<\infty \right\}.
\end{align*}
Notice that, as soon as $\nu>1$, if $\u\in\ell^1_\nu$ then the coefficients $u_k$ decay geometrically, and therefore the associated function $u$ is smooth.

The product of two functions $u$ and $v$ in function space gives rise to the discrete convolution product $\ast$ in sequence space:
\begin{align*}
\left(\u \ast \v\right)_n = \sum_{k\in\Z} u_{\vert k\vert} v_{\vert n-k\vert}.
\end{align*}
We recall that $\ell^1_\nu$ is a Banach algebra under the convolution product: for all $\u,\v$ in $\ell^1_\nu$,
\begin{align*}
\Vert \u\ast\v \Vert_\nu \leq \Vert \u \Vert_\nu \Vert \v \Vert_\nu.
\end{align*}
In particular, since $\Phi$ is assumed to be a polynomial, for any $\u$ in $\ell^1_\nu$ one can readily define $\Phi(\u)$, which also belongs to $\ell^1_\nu$, and similarly for $R$.
\begin{remark}
\label{rem:NP} 
In order for $R(\u)$ (or $\Phi(\u)$) to be well defined in $\ell^1_\nu$, it is sufficient to assume that $R$ is analytic on a disk of radius larger than $\Vert \u\Vert_\nu$. However, in some cases one might avoid the technicalities associated with having non polynomial terms by introducing extra variables and using automatic differentiation techniques, see e.g.~\cite{LesMirRan16}. A new alternative approach is also presented in Section~\ref{sec:NP}.
\end{remark} 

Given $\u\in\ell^1_\nu$, we denote by $M(\u):\ell^1_\nu \to \ell^1_\nu$ the associated multiplication operator, i.e. for all $v$ in $\ell^1_\nu$, $M(\u)v = u\ast v$.

Let $L$ be a linear operator on $\ell^1_\nu$. $L$ can be identified with an \emph{infinite dimensional matrix} $\left(L_{k,n}\right)_{k,n\geq 0}$ where, for any $\u$ in $\ell^1_\nu$,
\begin{align*}
(L\u)_k = \sum_{n=0}^\infty L_{k,n} u_n \qquad \forall~k\geq 0.
\end{align*}
We recall that the operator norm of $L$ can be easily expressed in terms of the coefficients $L_{k,n}$, namely
\begin{align*}
\Vert L\Vert_\nu := \sup_{\substack{\u\in\ell^1_\nu \\ \u\neq 0}} \frac{\Vert L\u\Vert_\nu}{\Vert \u\Vert_\nu} = \sup_{n\geq 0} \frac{1}{\xi_n(\nu)} \sum_{k=0}^\infty \vert L_{k,n}\vert \xi_k(\nu). 
\end{align*}
We point out that, for any $\u\in\ell^1_\nu$, the operator norm of the multiplication operator $M(\u)$ is nothing but the norm of the sequence $\u$, i.e.  $\Vert M(\u)\Vert_\nu = \Vert \u\Vert_\nu$. We also recall that the $\Vert\cdot\Vert_\nu$ norm of the coefficients controls the $\Vert\cdot\Vert_{C^0}$ norm of the associated function. That is, for any $\nu\geq 1$ and any $\u$ in $\ell^1_\nu$,
\begin{align}
\label{eq:C0l1}
\Vert u \Vert_{C^0} := \sup_{x\in\R} \vert u(x)\vert \leq \Vert \u \Vert_\nu.
\end{align}

We slightly abuse the notation $\Delta$ and also use it on a sequence of Fourier coefficients. That is, $\Delta\u$ is the sequence defined by
\begin{align*}
\left(\Delta \u\right)_n = -(n\pi)^2 u_n \qquad n\geq 0.
\end{align*}
Similarly $\Delta^{-1}\u$ is the sequence defined by
\begin{align*}
\left(\Delta^{-1} \u\right)_n = 
\left\{
\begin{aligned}
&0 \qquad & n=0, \\
&-\frac{1}{(n\pi)^2} u_n \qquad & n\geq 1.
\end{aligned}
\right.
\end{align*}
\begin{remark}
The Laplacian with Neumann boundary conditions is only invertible if we restrict ourselves to functions having zero mean, or in terms of Fourier coefficients if we only deal with modes $n\geq 1$. This will be enough for our purposes, and justifies the notation $\Delta^{-1}$ above.
\end{remark}

Finally, given $N\geq 1$ we introduce the projection onto the $N$ first modes $\Pi_N$, where
\begin{align*}
\left(\Pi_N \u \right)_n = 
\left\{
\begin{aligned}
&u_n \qquad & 0\leq n < N, \\
&0 \qquad & n\geq N,
\end{aligned}
\right.
\end{align*}
and the associated subspace $\Pi_N\ell^1_\nu$ of $\ell^1_\nu$. In the sequel, we frequently identify a vector $\bu=(u_0,\ldots,u_{N-1})$ of $\R^N$ and its injection $\bu = (u_0,\ldots,u_{N-1},0,\ldots,0,\ldots)$ in $\Pi_N\ell^1_\nu$. 

\subsection{$F=0$ into fixed point problem}
\label{sec:fixed-point}

Rewriting~\eqref{eq:PM} in Fourier space, our aim is to find $\u$ in $\ell^1_\nu$ satisfying
\begin{align*}
F_n(\u) := -(n\pi)^2 \Phi_n(\u) + R_n(\u) = 0 \qquad \forall~n\geq 0,
\end{align*}
or, in a more compact form
\begin{align}
\label{eq:F_PM}
F(\u) := \Delta \Phi(\u) + R(\u) = 0.
\end{align}
We now assume that we have an approximate zero $\bu\in\Pi_N\ell^1_\nu$ of $F$, which in practice can be obtained by numerically solving the finite dimensional problem $\Pi_N F\Pi_N = 0$.

In order to define a suitable approximate inverse $A$ of $DF(\bu)$, we first consider $\bw\in\Pi_N\ell^1_\nu$ such that $\bw \ast \Phi'(\bu) \approx \1$.
\begin{remark}
In practice, such $\bw$ is easily obtained by numerically solving the linear system 
\begin{align*}
\Pi_N M\left(\Phi'(\bu)\right) \Pi_N \bw = \1,
\end{align*}
identifying $\Pi_N M\left(\Phi'(\bu)\right) \Pi_N$ with an $N\times N$ matrix, and $\bw$ and $\1$ with vectors in $\R^N$.
\end{remark}

Next, we consider an $N\times N$ matrix $\bar A$, a numerically computed approximation of $\Pi_N \left ( DF(\bu)^{-1}\right) \Pi_N $, which we identify with an operator on $\Pi_N\ell^1_\nu$, and define the crucial operator $A$ as follows.
\begin{align}
\label{eq:A_PM}
A= \bar A + \left(M(\bw) \Delta^{-1} - \Pi_N \left( M(\bw) \Delta^{-1} \right) \Pi_N \right).
\end{align}
In order to better understand this operator, one might think of it as an \emph{infinite matrix}, which writes
\begin{align*}
\renewcommand\arraystretch{3} 
\newcommand{\ph}{\phantom{\frac{-1}{(N\pi)^2} \bw}}
\left(\begin{array}{ccccccccc}
 & & & & & \multicolumn{1}{|c}{0} & 0 & 0 & \ldots \\ 
 & & & & & \multicolumn{1}{|c}{\frac{-\bw_{N-1}}{(N\pi)^2} } & 0 & 0 & \ph \\ 
 & & \bar A & & & \multicolumn{1}{|c}{ \vdots} & \frac{-\bw_{N-1}}{((N+1)\pi)^2}  & 0 & \ddots \\ 
 & & & & & \multicolumn{1}{|c}{\frac{-\bw_2}{(N\pi)^2} } & \vdots & \frac{-\bw_{N-1}}{((N+2)\pi)^2}  & \ddots\\ 
\ph & \ph & \ph & \ph & \ph & \multicolumn{1}{|c}{\frac{-\bw_1}{(N\pi)^2} } & \frac{-\bw_2}{((N+1)\pi)^2}   & \vdots & \ddots \\ \cline{1-5}
0 & \frac{-\bw_{N-1}}{\pi^2} & \ldots & \frac{-\bw_2}{((N-2)\pi)^2} & \frac{-\bw_1}{((N-1)\pi)^2} & \frac{-\bw_0}{(N\pi)^2}  & \frac{-\bw_1}{((N+1)\pi)^2}   & \frac{-\bw_2}{((N+2)\pi)^2}  & \\
0 & 0 & \frac{-\bw_{N-1}}{(2\pi)^2} &\ldots & \frac{-\bw_2}{((N-1)\pi)^2} & \frac{-\bw_1}{(N\pi)^2} & \frac{-\bw_0}{((N+1)\pi)^2}  & \frac{-\bw_1}{((N+1)\pi)^2}  &  \ddots \\
0 & 0 & 0 & \frac{-\bw_{N-1}}{(3\pi)^2} &\ldots & \frac{-\bw_2}{(N\pi)^2} & \frac{-\bw_1}{((N+1)\pi)^2} & \frac{-\bw_0}{((N+2)\pi)^2}  &  \ddots\\
\vdots &  & \ddots & \ddots & \ddots & & \ddots &  \ddots & \ddots\\
\end{array}\right),
\end{align*}
or in a more compact form,
\begin{align*}
A= 
\left(\begin{array}{cccccc}
 & & & \multicolumn{1}{|c}{ } &  & \\ 
 & \bar A & & \multicolumn{1}{|c}{ } &  & \\ 
\phantom{\bar A} & \phantom{\bar A} & \phantom{\bar A} & \multicolumn{1}{|c}{ } &  & \\ \cline{1-3}
 & & & & & \\
 & & & & M(\bw)\Delta^{-1} & \\
 & & & & &
\end{array}\right).
\end{align*}

\begin{remark}
\label{rem:barA}
When there is only linear diffusion, $\bar A$ is usually defined by numerically computing the inverse of $\Pi_N DF(\bu)\Pi_N$. In our situation this won't necessarily be good enough. This is related to the fact that the extra-diagonal terms in $DF(\bu)$ cannot be neglected, even outside of the finite block corresponding to the projection on $\Pi_N\ell^1_\nu$. For more details and an explicit computation, see Appendix~\ref{sec:barA}. 

Instead, in order to get a good enough approximation of $\Pi_N \left ( DF(\bu)^{-1}\right) \Pi_N $, we numerically compute the inverse of a larger block, and then project it back onto $\Pi_N\ell^1_\nu$, e.g.
\begin{align*}
\bar A \approx \Pi_N \left(\Pi_{2N} DF(\bu)\Pi_{2N}\right)^{-1} \Pi_N, 
\end{align*}
where the $\approx$ sign is only here to emphasize that the inverse does not have to be computed exactly.
\end{remark}

Now that the important operator $A$ has been defined, we are ready to consider the fixed-point operator $T:\ell^1_\nu\to\ell^1_\nu$ defined by
\begin{align}
\label{eq:T_PM}
T(u) = u - AF(u),
\end{align}
and to give the classical sufficient conditions for $T$ to admit a fixed-point near $\bu$. This Newton-Kantorovich-like theorem tells us that $A$ has to be a \emph{sufficiently good} approximate inverse of $DF(\bu)$, and that $\bu$ has to be a \emph{sufficiently good} approximate zero of $F$, makes it precise what \emph{sufficiently good} means (see~\eqref{eq:cond}), and then gives an explicit error bound for $\bu$ (see~\eqref{eq:r}).

\begin{theorem}
\label{th:NK}
With the notations introduced in this section, assume there exist constants $Y$, $Z_1$ and $Z_2$ satisfying
\begin{subequations}
\label{eq:YZ12}
\begin{align}
\Vert AF(\bu) \Vert_\nu &\leq Y \label{eq:Y}\\
\Vert I - ADF(\bu) \Vert_\nu  &\leq Z_1 \label{eq:Z1}\\
\Vert AD^2F(\u) \Vert_\nu &\leq Z_2 \qquad \forall~\u\in\ell^1_\nu, \label{eq:Z2}
\end{align}
\end{subequations}
and
\begin{subequations}
\label{eq:cond}
\begin{align}
Z_1 &< 1 \label{eq:Z1l1}\\
2YZ_2 &< (1-Z_1)^2.
\end{align}
\end{subequations}
Then, for any $r$ satisfying 
\begin{align}
\label{eq:r}
\frac{1-Z_1 - \sqrt{(1-Z_1)^2-2YZ_2}}{Z_2} \leq r < \frac{1-Z_1}{Z_2},
\end{align}
there exists a unique fixed-point $\u^*$ of $T$ in $\B_\nu(\bu,r)$, the closed ball of center $\bu$ and radius $r$ in $\ell^1_\nu$.

Assume further that $\bw$ (which plays a role in the definition of $A$), is such that
\begin{align}
\label{eq:hyp_w}
\Vert \1 - \bw\ast\Phi'(\bu) \Vert_\nu <1.
\end{align}
Then $\u^*$ is the unique zero of $F$ in $\B_\nu(\bu,r)$.
\end{theorem}
\begin{proof}
The first part of the proof, which consists in showing that $T$ defined in~\eqref{eq:T_PM} admits a unique fixed point in $\B_\nu(\bu,r)$ is standard, so we only sketch it. For any $\u$ in $\B_\nu(\bu,r)$ and $r$ satisfying~\eqref{eq:r}, we have
\begin{align*}
\Vert T(\u) - \bu \Vert_\nu \leq Y + Z_1 r + \frac{1}{2}Z_2 r^2 \leq r,
\end{align*}
and 
\begin{align*}
\Vert DT(\u) \Vert_\nu \leq Z_1 + Z_2 r <1,
\end{align*}
therefore $T$ is a contraction on $\B_\nu(\bu,r)$, and there is a unique fixed-point $\u^*$ of $T$ in $\B_\nu(\bu,r)$. 

The second part of the proof consists in showing that $A$ is injective, so that this fixed point of $T$ indeed corresponds to a zero of $F$. From~\eqref{eq:Z1} and~\eqref{eq:Z1l1}, we know that $ADF(\bu)$ is invertible, therefore we have at least that $A$ is surjective. Furthermore, by~\eqref{eq:hyp_w} we have that 
\begin{align*}
\Vert I - M(\bw)\, M(\Phi'(\bu)) \Vert_\nu = \Vert \1 - \bw\ast\Phi'(\bu) \Vert_\nu <1,
\end{align*}
therefore $M(\bw)$ is invertible ($M(\bw)$ and $M(\Phi'(\bu))$ commute, because $\ast$ is associative). Since $\Delta^{-1}$ is a Fredholm operator of index $0$, so is $M(\bw) \Delta^{-1}$. However, $A$ is just a compact perturbation of $M(\bw) \Delta^{-1}$, therefore it is also Fredholm of index $0$, and since we have already shown that $A$ is surjective, it must be injective.
\end{proof}

\begin{remark}
\label{rem:thNK}
With a Fourier discretization and linear diffusion, the second part of the above proof is usually trivial, but ideas similar to the ones used here already appeared in the context of computer-assisted proofs, see e.g. the proof of~\cite[Theorem 6.2]{NakPluWat19}. 

While this second part of the proof seems to require an extra assumption, namely~\eqref{eq:hyp_w}, the term $\Vert \1 - \bw\ast\Phi'(\bu) \Vert_\nu$ is going to naturally appear in our estimate for $Z_1$, and therefore~\eqref{eq:hyp_w} will be automatically satisfied as soon as~\eqref{eq:Z1l1} is, see~\eqref{eq:Z1tail_PM}.

Finally let us mention that~\eqref{eq:Z2} does not really need to hold for all $\u$ in $\ell^1_\nu$, but only in a small neighborhood of $\bu$. In the present case it does not matter, since we took $\Phi$ and $R$ as polynomials of degree two and therefore $D^2F$ is constant. In a more general case one may introduce an a priori radius $r^*$, and restrict the whole analysis to $\B_\nu(\bu,r^*)$, in order to only have to estimate $D^2F(\u)$ locally.
\end{remark}

\subsection{Derivation of the bounds}
\label{sec:boundsPM}

It remains to derive computable estimates $Y$, $Z_1$ and $Z_2$ satisfying~\eqref{eq:YZ12}, and to check that~\eqref{eq:cond} also holds. 

One key point to keep in mind in this subsection is that, since $\bw\in\Pi_N\ell^1_\nu$, the infinite matrix $M(\bw)$ is a banded matrix, with a bandwidth of at most $N$. That is, as soon as $\vert k-n\vert \geq N$, $\left(M(\bw)\right)_{k,n}=0$. Since we assumed $\Phi(u)=u^2$ and $R(u) = \alpha u -\beta u^2 + g$, the same is true for $M(\Phi'(\bu))$ and $M(R'(\bu))$.

\begin{remark}
If $R$ or $\Phi$ were polynomials of higher order, $M(\Phi'(\bu))$ and $M(R'(\bu))$ would simply have a larger bandwidth (equal to $d(N-1)+1$ if the polynomial is of degree $d$), and all the estimates to come would have to be adapted in a straightforward manner.

If $R$ (or $\Phi$) were to be merely analytic, one could split $R'(\bu)$ between a finite sequence and a remainder whose $\Vert \cdot \Vert_\nu$ norm is small and could be estimated, deal as above with the finite part, and keep track of the extra terms that are produced by the remainder. See also Section~\ref{sec:NP} for a different approach.

Still for the sake of simplicity, in the sequel we also assume that $\g\in\Pi_N\ell^1_\nu$, but more general $g$ could be handled in a similar way, as soon as the remainder can be estimated explicitly.
\end{remark}

\subsubsection{The bound $Y$}
\label{sec:Y_PM}

Since $\Phi$ and $R$ are polynomials, and $\bu$ and $\g$ belong to $\Pi_N\ell^1_\nu$, $F(\bu)$ only has a finite number of non-zero coefficients. Similarly, since $\bw$ belongs to $\Pi_N\ell^1_\nu$, each column of $A$ only has a finite number of non-zero coefficients. Therefore $AF(\u)$ and then $\Vert AF(\u) \Vert_\nu$ can be computed exactly, up to rounding errors, and the rounding errors can be controlled using interval arithmetic. In our implementation we make use of the Intlab package~\cite{Rum99} for Matlab. The upper bound of this computation will be our bound $Y$, and satisfies~\eqref{eq:Y}.

\subsubsection{The bound $Z_1$}
\label{sec:Z1_PM}

Denoting $B=I - ADF(\bu)$ and $B_{(\cdot,n)}$ the $n$-th column of $B$, we split the estimation of the norm of $B$ as follows.
\begin{align*}
\Vert B\Vert_\nu &= \sup_{n\geq 0} \frac{1}{\xi_n(\nu)} \Vert B_{(\cdot,n)} \Vert_\nu \\
&= \max\left( \max_{0\leq n \leq 2N-2} \frac{1}{\xi_n(\nu)} \Vert B_{(\cdot,n)}\Vert_\nu,\ \sup_{n\geq 2N-1} \frac{1}{\xi_n(\nu)} \Vert B_{(\cdot,n)} \Vert_\nu \right) \\
&\leq \max\left( Z_1^{finite},\ Z_1^{tail} \right).
\end{align*}
$Z_1^{finite}$ can be computed explicitly. Indeed, the $n$-th column of $B$ is nothing but the $n$-th column of $DF(\bu)$ multiplied by $A$, and similarly to the situation for the $Y$ bounds, the norm of any individual column can be computed exactly (or, to be more precise, enclosed rigorously using interval arithmetic). 

Regarding the estimate for $Z_1^{tail}$, we use the fact that $M(\Phi'(\bu))$ and $M(R'(\bu))$ have bandwidth $N$. Hence for any column of $DF(\bu)=\Delta M(\Phi'(\bu)) + M(R'(\bu))$ having index $n\geq 2N-1$, its $N$ first coefficients must be zero, and the multiplication of such a column by $A$ does not depend on $\bar A$. This enables us to easily estimate by hand $\Vert B_{(\cdot,n)} \Vert_\nu$ for any $n\geq 2N-1$. Equivalently, we can start by rewriting
\begin{align*}
\sup_{n\geq 2N-1} \frac{1}{\xi_n(\nu)} \Vert B_{(\cdot,n)} \Vert_\nu = \sup\limits_{\substack{\u\in(I-\Pi_{2N-1})\ell^1_\nu \\ \u\neq 0 }} \frac{\Vert B \u \Vert}{\Vert \u\Vert}.
\end{align*}
Then, since $\Phi'(\bu)$ and $R'(\bu)$ belong to $\Pi_N\ell^1_\nu$, for any $\u\in(I-\Pi_{2N-1})\ell^1_\nu$ we have that $\Phi'(\bu)\ast \u$ and $R'(\bu)\ast \u$ belong to $(I-\Pi_N)\ell^1_\nu$, and so does $DF(\bu)\u$. Therefore, for any $\u\in(I-\Pi_{2N-1})\ell^1_\nu$,
\begin{align*}
B\u &= \u - ADF(\bu)\u \\
&= \u -\bw\ast \Delta^{-1}DF(\bu)\u \\
&= \u -\bw\ast \Phi'(\bu)\ast \u - \bw\ast\Delta^{-1}\left(R'(\bu)\ast\u\right),
\end{align*}
and 
\begin{align*}
\Vert B\u\Vert_\nu &\leq  \Vert \1 -\bw\ast \Phi'(\bu)\Vert_\nu \Vert \u\Vert_\nu + \Vert \bw\Vert_\nu \Vert\Delta^{-1}\left(R'(\bu)\ast\u\right)\Vert_\nu.
\end{align*}
Finally, since $R'(\bu)\ast\u$ belongs to $(I-\Pi_N)\ell^1_\nu$ we have $\Vert\Delta^{-1}\left(R'(\bu)\ast\u\right)\Vert_\nu \leq \frac{1}{(N\pi)^2} \Vert R'(\bu)\ast\u\Vert_\nu$, and we can take
\begin{align}
\label{eq:Z1tail_PM}
Z_1^{tail} = \Vert \1 - \bw\ast\Phi'(\bu)\Vert_\nu + \frac{1}{(N\pi)^2} \Vert \bw\ast R'(\bu)\Vert_\nu.
\end{align}
We then just take $Z_1 = \max\left( Z_1^{finite},\ Z_1^{tail} \right)$, and it satisfies~\eqref{eq:Z1}.

\subsubsection{The bound $Z_2$}
\label{sec:Z2_PM}

For any $\u_1,\u_2$ in $\ell^1_\nu$, remembering that $\Phi(u)=u^2$ and $R(u) = \alpha u - \beta u^2 + g$, we have
\begin{align*}
D^2F(\u)(\u_1,\u_2) = 2\Delta \left(\u_1\ast\u_2\right) -2\beta \u_1\ast\u_2,
\end{align*}
and therefore
\begin{align}
\label{eq:AD2F}
\Vert AD^2F(\u) \Vert_\nu &= \sup_{\substack{\u_1,\u_2\in\ell^1_\nu \\ \u_1,\u_2\neq 0}} \frac{\Vert AD^2F(\u)(\u_1,\u_2) \Vert_\nu}{\Vert \u_1\Vert_\nu \Vert \u_2\Vert_\nu} \nonumber \\
&\leq  2\Vert A\Delta \Vert_\nu + 2\vert\beta\vert \Vert A\Vert_\nu .
\end{align}
We thus take
\begin{align*}
Z_2 =   2\Vert A\Delta \Vert_\nu + 2\vert\beta\vert \Vert A\Vert_\nu .
\end{align*}

Similarly to what we have done for the $Z_1$ bound, the norm of $A$ and $A\Delta$ can be computed by splitting between an finite number of column and the rest, and we get
\begin{subequations}
\label{eq:norm_A}
\begin{align}
\Vert A\Vert_\nu &= \max\left( \max_{0\leq n \leq N-1} \frac{1}{\xi_n(\nu)} \Vert A_{(\cdot,n)}\Vert_\nu, \frac{1}{(N\pi)^2}\Vert \bw \Vert_\nu \right) \\
\Vert A\Delta\Vert_\nu &= \max\left( \max_{0\leq n \leq N-1} \frac{1}{\xi_n(\nu)} \Vert \left(A\Delta\right)_{(\cdot,n)}\Vert_\nu, \Vert \bw \Vert_\nu \right).
\end{align}
\end{subequations}
Combining~\eqref{eq:AD2F} and~\eqref{eq:norm_A} gives us a computable $Z_2$ bound satisfying~\eqref{eq:Z2}.

\subsection{Example and results}
\label{sec:res_PM}

\begin{figure}
\centering
\begin{subfigure}[t]{0.49\textwidth}
\centering
\includegraphics[scale=0.49]{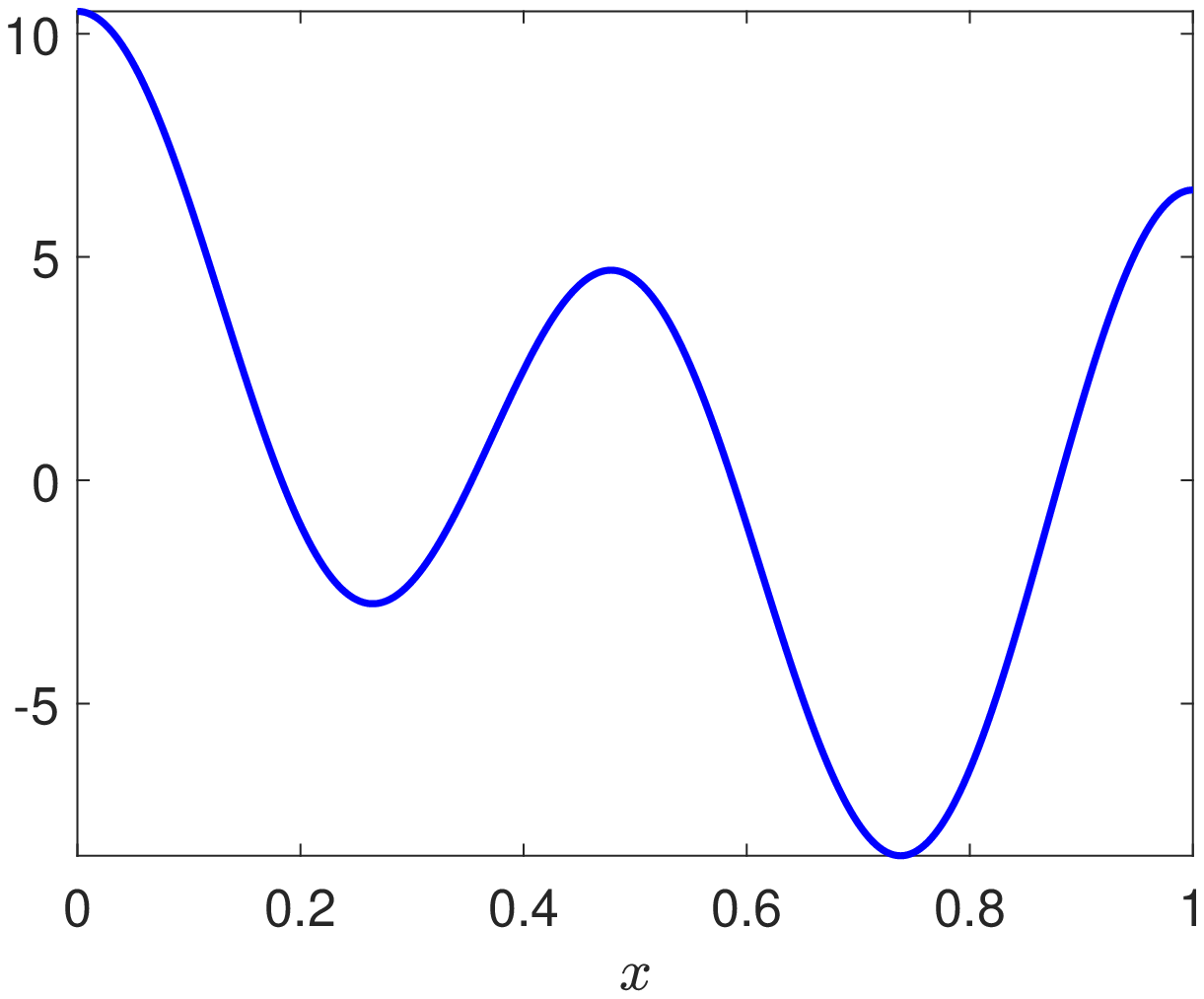}
\caption{The forcing term $g$ used in the examples.}
\label{fig:g}
\end{subfigure}
\begin{subfigure}[t]{0.49\textwidth}
\centering
\includegraphics[scale=0.49]{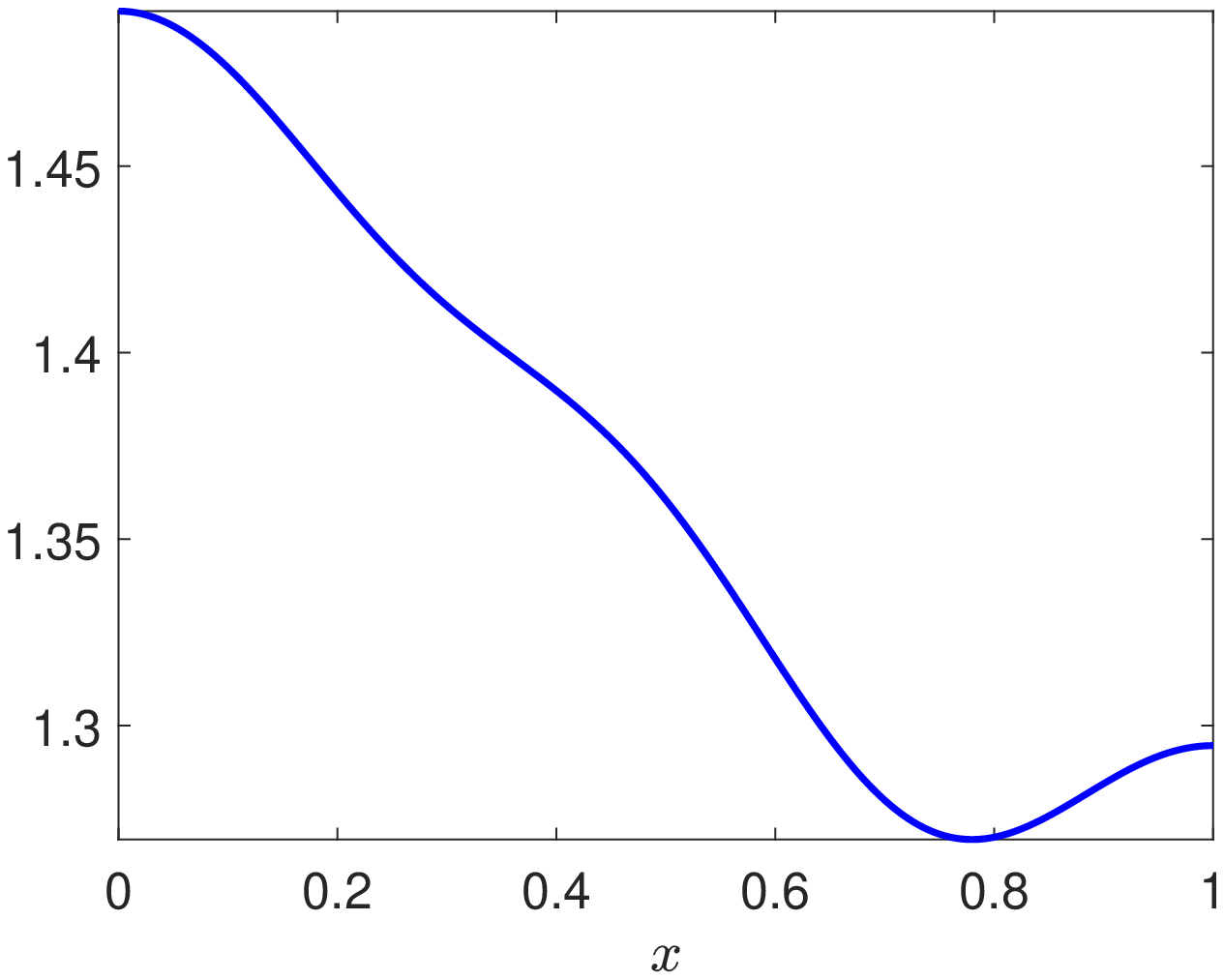}
\caption{The approximate solution $\bar u$ of~\eqref{eq:PM} with $\Phi(u) = u^2$, $R(u) = u-u^2 +g$ and $g$ as in Figure~\ref{fig:g}, which has been validated in Theorem~\ref{th:PM}.}
\label{fig:u_PM}
\end{subfigure}
\caption{Data associated to Theorem~\ref{th:PM}.}
\end{figure}

In this subsection, we consider an explicit case, namely $\Phi(u)=u^2$ and $R(u) = u - u^2 + g$ (i.e. we take $\alpha=\beta=1$), and $g$ given by $g(x) = \frac{1}{2} + 3\cos(\pi x) +2\cos(2\pi x) - \cos(3\pi x) + 6\cos(4\pi x)$ (see Figure~\ref{fig:g}), for which we check that Theorem~\ref{th:NK} is applicable.

We computed an approximate solution $\bu$ (see Figure~\ref{fig:u_PM}), and then validated it using the procedure described in this whole section, which yield the following theorem.

\begin{theorem}
\label{th:PM}
Let $\bar u $ be the function whose Fourier coefficients $\bu$ can be downloaded at~\cite{Bre21}, and which is represented in Figure~\ref{fig:u_PM}. There exists a strong solution $u$ of~\eqref{eq:PM}, with $\Phi$ and $R$ as given just above, such that $\left\Vert u-\bar u\right\Vert_{C^0} \leq 2\times 10^{-10}$.
\end{theorem}
\begin{proof}
We apply the validation procedure with $N=20$ and $\nu=1.1$. That is, we define some $\bw$ and $\bar A$, then evaluate the bounds $Y$, $Z_1$ and $Z_2$ obtained in Section~\ref{sec:boundsPM}, with interval arithmetic to take rounding errors into accounts. We get that assumption~\eqref{eq:YZ12} is satisfied with
\begin{align*}
Y= 1.3\times 10^{-10},\quad Z_1 = 0.0002 \quad \text{and} \quad Z_2 = 2.2.
\end{align*}
Therefore~\eqref{eq:cond} holds, and so does~\eqref{eq:r} with $r=2\times 10^{-10}$. As mentioned in Remark~\ref{rem:thNK}, formula~\eqref{eq:Z1tail_PM} show that~\eqref{eq:hyp_w} also holds because~\eqref{eq:Z1l1} does. 

We can then apply Theorem~\ref{th:NK}, and there exists a unique zero $\u$ of $F$ in $\B_\nu(\bu,r)$. Since $\nu>1$ the associated function $u$ is smooth, and is therefore a strong solution of~\eqref{eq:PM}. The announced error estimate in $C^0$ norm follows directly from~\eqref{eq:C0l1}.

The computational parts of the proof can be reproduced by running the Matlab code \texttt{script\_PM.m} from~\cite{Bre21}, together with Intlab~\cite{Rum99}.
\end{proof}
\begin{remark}
\label{rem:space}
If needed, one could also get error estimates in different function spaces, since the $\ell^1_\nu$ norm of $\u$ controls the norm of any derivative of $u$, as soon as $\nu>1$ (with constants depending on $\nu$).

In $\ell^1_\nu$ we not only have existence of a zero $\u$ of $F$, but also local uniqueness. This local uniqueness carries over to the solutions of~\eqref{eq:PM}, provided we know a priori that any solution of~\eqref{eq:PM} must be smooth enough for its Fourier coefficients to belong to $\ell^1_\nu$. Otherwise, we only have local uniqueness among smooth enough functions.

For problems where analyticity is too much too ask for, be it because the solution itself is not analytic, or because we care about local uniqueness but proving the analyticity a priori is hard, it should be noted that one can easily replicated the whole procedure presented here in a sequence space corresponding to functions of lower regularity, see e.g.~\cite{LesMir17}.
\end{remark}

\section{The SKT system}
\label{sec:SKT}

We now use the techniques introduced in Section~\ref{sec:PM} to study steady states of the SKT system
\begin{align}
\label{eq:SKT}
\left\{
\begin{aligned}
&\partial_t u^{(1)} = \Delta \left(\left(d_1 + d_{11} u^{(1)} + d_{12} u^{(2)} \right) u^{(1)}\right)  + \left(r_1 - a_1 u^{(1)} - b_1 u^{(2)} \right) u^{(1)}  = 0 \qquad & \text{on }\Omega,  \\
&\partial_t u^{(2)} = \Delta \left(\left(d_2 + d_{21} u^{(1)} + d_{22} u^{(2)} \right) u^{(2)}\right)  + \left(r_2 - b_2 u^{(1)} - a_2 u^{(2)} \right) u^{(2)}  = 0 \qquad & \text{on }\Omega,  \\
&\frac{\partial u^{(1)}}{\partial n } = 0 = \frac{\partial u^{(2)}}{\partial n } \qquad & \text{on }\partial\Omega,
\end{aligned}
\right.
\end{align}
which is exactly of the form~\eqref{eq:NLD}, with
\begin{align*}
u= \begin{pmatrix}
u^{(1)} \\ u^{(2)}
\end{pmatrix} ,
\quad
\Phi(u) = \begin{pmatrix}
\Phi^{(1)}(u) \\ \Phi^{(2)}(u)
\end{pmatrix} =
\begin{pmatrix}
\left(d_1 + d_{11} u^{(1)} + d_{12} u^{(2)} \right) u^{(1)} \\
\left(d_2 + d_{21} u^{(1)} + d_{22} u^{(2)} \right) u^{(2)}
\end{pmatrix},
\end{align*}
and
\begin{align*}
R(u) = \begin{pmatrix}
R^{(1)}(u) \\ R^{(2)}(u)
\end{pmatrix} =
\begin{pmatrix}
\left(r_1 + a_1 u^{(1)} + b_1 u^{(2)} \right) u^{(1)} \\
\left(r_2 + b_2 u^{(1)} + a_2 u^{(2)} \right) u^{(2)}
\end{pmatrix}.
\end{align*}

\subsection{Presentation of the model}
\label{sec:SKTintro}

The SKT system was introduced in the seminal paper~\cite{ShiKawTer79} in order to model the dynamics of two competing species, whose densities are denoted here by $u^{(1)}$ and $u^{(2)}$. The reactions terms are standard Lotka-Volterra terms with signs indicating intra-specific and inter-specific competitions. The key feature of this model is the presence of nonlinear diffusion terms: the diffusion rate of each species depends on the local density of both species. These nonlinear diffusion terms can give rise to a repulsive effect which leads to \emph{spatial segregation}: the two species co-exist but mostly concentrate in different regions of the domain $\Omega$. Mathematically, this corresponds to the existence of non-homogeneous steady states of~\eqref{eq:SKT} exhibiting this type of patterns.

The existence and stability analysis of such steady states has been studied extensively, through a wide variety of techniques such as bifurcation theory, singular perturbation theory or fixed point index theory, see e.g.~\cite{MimKaw80,MimNisTesTsu84,LouNi96,RyuAhn03}, which provide a qualitative understanding on the conditions required on the many parameters of~\eqref{eq:SKT} for non-homogeneous steady states to exist. Some asymptotic parameter regimes have also been scrutinized, especially when one or both of the cross-diffusion coefficients $d_{12}$ and $d_{12}$ become very large, and in such cases additional information about the shape of the obtained solutions is also available, see e.g.~\cite{Ni98,LouNiYot04}. However, numerical studies like~\cite{IidMimNim06,IzuMim08} show that the steady states of~\eqref{eq:SKT} can be very diverse, already in the one dimensional case, and that many of them can co-exists for a give set of parameter values. Proving the existence, and characterizing the shape, of all these different steady states seems very hard by purely analytic means, but computer-assisted techniques provide a valuable tool to attack such questions.

The computer-assisted study of the steady states of the SKT system started first via a system with linear diffusion approximating~\eqref{eq:SKT} in~\cite{BreLesVan13}, and then for~\eqref{eq:SKT} itself  in~\cite{BreCas18}, but only in the so-called \emph{triangular case}. This case corresponds to taking the \emph{self-diffusion} coefficients $d_{11}$ and $d_{22}$ equal to $0$, and one of the cross-diffusion coefficient (say $d_{21}$) equal to $0$. While having only $d_{12}uv$ as a nonlinear diffusion term is already sufficient to produce many interesting steady states, the mathematical analysis is somewhat simpler in that case. In particular, the work~\cite{BreCas18} takes advantage of the fact that, in the triangular case, $\Phi^{-1}$ can be computed explicitly, and therefore one can directly transform~\eqref{eq:SKT} into a system with only linear diffusion. However, this feature is specific to the triangular case, and the general case where $d_1$, $d_2$, $d_{12}$, $d_{21}$, $d_{11}$ and $d_{22}$ all are possibly non-zero could not be handled in~\cite{BreCas18}. 

With the techniques introduced in this paper, we get a new approach which is not only more general (i.e. not restricted to the triangular case), but also more efficient compared to the ad hoc technique used in~\cite{BreCas18} for the triangular case (see Remark~\ref{rem:1_N2_SKT}). Building upon the work done in Section~\ref{sec:PM}, we introduce the necessary supplementary notations and spaces in Section~\ref{sec:not_SKT}, the fixed-point setup in Section~\ref{sec:fp_SKT}, derive the estimates needed for the validation in Section~\ref{sec:bounds_SKT}, and present some examples in Section~\ref{sec:res_SKT}, which include the a posteriori validation of some interesting solutions first found numerically in~\cite{BreKueSor21}.

\subsection{Notations}
\label{sec:not_SKT}

We consider the product space $\XX_\nu:=\ell^1_\nu\times\ell^1_\nu$, and slightly abuse all the notations introduced in Section~\ref{sec:not} by applying them component-wise. That is, 
\begin{align*}
\renewcommand{\arraystretch}{1.5}
\text{for any }
\u = \begin{pmatrix}
\u^{(1)} \\ \u^{(2)}
\end{pmatrix},\quad
\Pi_N\u := \begin{pmatrix}
\Pi_N\u^{(1)} \\ \Pi_N\u^{(2)}
\end{pmatrix}, \quad 
\Vert \u \Vert_\nu := \begin{pmatrix}
\Vert \u^{(1)}\Vert_\nu \\ \Vert \u^{(2)}\Vert_\nu
\end{pmatrix}, \quad 
M(\u) := \begin{pmatrix}
M(\u^{(1)}) \\ M(\u^{(2)})
\end{pmatrix}.
\end{align*}
The norm on $\XX_\nu$ is then defined as
\begin{align*}
\left\Vert \u \right\Vert_{\XX_\nu} = \left\vert \Vert \u \Vert_\nu\right\vert_1 = \Vert\u^{(1)}\Vert_\nu + \Vert\u^{(2)}\Vert_\nu,
\end{align*}
where $\left\vert\cdot\right\vert_1$ denotes the $1$-norm on $\R^2$.

Given a linear operator $L$ on $\XX_\nu$, written
\begin{align*}
L = 
\renewcommand{\arraystretch}{2}
\left(
\begin{array}{c|c}
L^{(1,1)} & L^{(1,2)} \\ \hline
L^{(2,1)} & L^{(2,2)}
\end{array}
\right),
\end{align*}
where each $L^{(i,j)}$ is an operator on $\ell^1_\nu$, we also abuse $\Vert\cdot\Vert_\nu$ notation in a similar way, that is
\begin{align*}
\renewcommand{\arraystretch}{1.5}
\Vert L \Vert_\nu = \begin{pmatrix}
\Vert L^{(1,1)}\Vert_\nu & \Vert L^{(1,2)}\Vert_\nu \\
\Vert L^{(2,1)}\Vert_\nu & \Vert L^{(2,2)}\Vert_\nu
\end{pmatrix}.
\end{align*}
For the operator norm  of $L$, it is then straightforward to show that
\begin{align*}
\Vert L\Vert_{\XX_\nu} &\leq \max\left(\Vert L^{(1,1)}\Vert_\nu + \Vert L^{(2,1)}\Vert_\nu ,\, \Vert L^{(1,2)}\Vert_\nu + \Vert L^{(2,2)}\Vert_\nu \right) \\
& = \left\vert \Vert L \Vert_\nu \right\vert_1.
\end{align*}

\subsection{$F=0$ into fixed point problem}
\label{sec:fp_SKT}

We look for a zero $\u=(\u^{(1)},\u^{(2)})$ in $\XX_\nu$ of
\begin{align*}
F=\begin{pmatrix}
F^{(1)} \\
F^{(2)}
\end{pmatrix} 
\end{align*}
defined as
\begin{align*}
\left\{
\begin{aligned}
F^{(1)}_n(\u) &= -(N\pi)^2 \Phi^{(1)}_n(\u) + R^{(1)}_n(\u) \\
F^{(2)}_n(\u) &= -(N\pi)^2 \Phi^{(2)}_n(\u) + R^{(2)}_n(\u) 
\end{aligned}
\right.
\qquad \forall~n\geq 0,
\end{align*}
or, in a more condensed form,
\begin{align*}
F(\u) = \Delta\Phi(\u) + R(\u),
\end{align*}
with $\Phi$ and $R$ as in Section~\ref{sec:SKTintro}.

We now assume that we have computed an approximate zero $\bu$ in $\Pi_N\XX_\nu$ that we want to validate a posteriori. We also assume that we have computed
\begin{align*}
\bw = \begin{pmatrix}
\bw^{(1,1)} & \bw^{(1,2)} \\ \bw^{(2,1)} & \bw^{(2,2)}
\end{pmatrix} \in \M_{2,2}(\Pi_N\ell^1_\nu), \quad \text{such that}\quad \bw \ast D\Phi (\bu) \approx \begin{pmatrix}
\1 & \0 \\ \0 & \1
\end{pmatrix},
\end{align*}
where
\begin{align*}
D\Phi(\bu) = \begin{pmatrix}
d_1 + 2d_{11} \bu^{(1)} + d_{12} \bu^{(2)} & d_{12} \bu^{(1)} \\
d_{21} \bu^{(2)} & d_2 + d_{21} \bu^{(1)} + 2d_{22} \bu^{(2)}
\end{pmatrix},
\end{align*}
and 
$\bw \ast D\Phi (\bu)$ must be understood as
\begin{align*}
\begin{small}
\begin{pmatrix}
\bw^{(1,1)} \ast \left(d_1 + 2d_{11} \bu^{(1)} + d_{12} \bu^{(2)}\right) + \bw^{(1,2)} \ast d_{21} \bu^{(2)} &  \bw^{(2,1)} \ast \left(d_1 + 2d_{11} \bu^{(1)} + d_{12} \bu^{(2)}\right) + \bw^{(2,2)} \ast d_{21} \bu^{(2)}\\
\bw^{(1,1)} \ast d_{12} \bu^{(1)} + \bw^{(1,2)} \ast \left(d_2 + d_{21} \bu^{(1)} + 2d_{22} \bu^{(2)}\right) & \bw^{(2,1)} \ast d_{12} \bu^{(1)} + \bw^{(2,2)} \ast \left(d_2 + d_{21} \bu^{(1)} + 2d_{22} \bu^{(2)}\right)
\end{pmatrix}.
\end{small}
\end{align*}
Next, we consider $\bar A$, a numerically computed approximation of $\Pi_N \left ( DF(\bu)^{-1}\right) \Pi_N $ (Remark~\ref{rem:barA} is also relevant here). Identifying $\bar A$ with a $2N\times 2N$ matrix, we separate it into four $N\times N$ blocks
\begin{align*}
\bar A = 
\renewcommand{\arraystretch}{2}
\left(
\begin{array}{c|c}
\bar A^{(1,1)} & \bar A^{(1,2)} \\ \hline
\bar A^{(2,1)} & \bar A^{(2,2)}
\end{array}
\right).
\end{align*}
The operator $A$ is then defined as
\begin{align*}
\newcommand{\phs}{\phantom{\bar A}}
\newcommand{\phl}{\phantom{M(\sigma)\Delta^{-1}}}
A = 
\left(
\begin{array}{c|c}
\begin{array}{cccccc}
 & & & \multicolumn{1}{|c}{ } &  & \\ 
 & \bar A^{(1,1)} & & \multicolumn{1}{|c}{ } &  & \\ 
\phs & \phs & \phs & \multicolumn{1}{|c}{ } &  & \\ \cline{1-3}
 & & & & & \\
 & & & & M\left(\bw^{(1,1)}\right)\,\Delta^{-1} & \\
 & & & & \phl &
\end{array}
&
\begin{array}{cccccc}
 & & & \multicolumn{1}{|c}{ } &  & \\ 
 & \bar A^{(1,2)} & & \multicolumn{1}{|c}{ } &  & \\ 
\phs & \phs & \phs & \multicolumn{1}{|c}{ } &  & \\ \cline{1-3}
 & & & & & \\
 & & & & M\left(\bw^{(1,2)}\right)\,\Delta^{-1} & \\
 & & & & &
\end{array} \\ \hline
\begin{array}{cccccc}
 & & & \multicolumn{1}{|c}{ } &  & \\ 
 & \bar A^{(2,1)} & & \multicolumn{1}{|c}{ } &  & \\ 
\phs & \phs & \phs & \multicolumn{1}{|c}{ } &  & \\ \cline{1-3}
 & & & & & \\
 & & & & M\left(\bw^{(2,1)}\right)\,\Delta^{-1} & \\
 & & & & \phl &
\end{array}
&
\begin{array}{cccccc}
 & & & \multicolumn{1}{|c}{ } &  & \\ 
 & \bar A^{(2,2)} & & \multicolumn{1}{|c}{ } &  & \\ 
\phs & \phs & \phs & \multicolumn{1}{|c}{ } &  & \\ \cline{1-3}
 & & & & & \\
 & & & & M\left(\bw^{(2,2)}\right)\,\Delta^{-1} & \\
 & & & & \phl &
\end{array}
\end{array}
\right),
\end{align*}
or, in a more compact form,
\begin{align}
A= \bar A + \left(M(\bw) \Delta^{-1} - \Pi_N \left( M(\bw) \Delta^{-1} \right) \Pi_N \right).
\end{align}

Now that $A$ is defined, we can again state sufficient conditions for the validation, which are the same as the one in Theorem~\ref{th:NK}, up to the slightly different space and norm.
\begin{theorem}
\label{th:NK_SKT}
With the notations introduced in this section, assume there exist constants $Y$, $Z_1$ and $Z_2$ satisfying
\begin{align*}
\Vert AF(\bu) \Vert_{\XX_\nu} &\leq Y \\
\Vert I - ADF(\bu) \Vert_{\XX_\nu}  &\leq Z_1 \\
\Vert AD^2F(\u) \Vert_{\XX_\nu} &\leq Z_2 \qquad \forall~\u\in\XX_\nu, 
\end{align*}
and
\begin{align*}
Z_1 &< 1 \\
2YZ_2 &< (1-Z_1)^2.
\end{align*}
Then, for any $r$ satisfying 
\begin{align*}
\frac{1-Z_1 - \sqrt{(1-Z_1)^2-2YZ_2}}{Z_2} \leq r < \frac{1-Z_1}{Z_2},
\end{align*}
there exists a unique fixed-point $\u^*$ of $T$ in $\B_{\XX_\nu}(\bu,r)$, the closed ball of center $\bu$ and radius $r$ in $\XX_\nu$.
Assume further that $\bw$ (which plays a role in the definition of $A$), is such that
\begin{align*}
\left\Vert \begin{pmatrix}
\1 & \0 \\ \0 & \1
\end{pmatrix} - \bw \ast D\Phi (\bu) \right\Vert_{\XX_\nu} <1.
\end{align*}
Then $\u^*$ is the unique zero of $F$ in $\B_{\XX_\nu}(\bu,r)$.
\end{theorem}

\subsection{Derivation of the bounds}
\label{sec:bounds_SKT}

It remains to derive computable estimates $Y$, $Z_1$ and $Z_2$ satisfying the assumptions of Theorem~\ref{th:NK_SKT}. Since all the estimations are very similar to the one performed in Section~\ref{sec:boundsPM}, we omit the details.

\subsubsection{The bound $Y$}

As in Section~\ref{sec:Y_PM}, there is no pen and paper estimation to make here, $\Vert AF(\bu)\Vert_{\XX_\nu}$ can be computed explicitly using interval arithmetic.

\subsubsection{The bound $Z_1$}

As in Section~\ref{sec:Z1_PM}, we introduce
\begin{align*}
B=\begin{pmatrix}
B^{(1,1)} & B^{(1,2)} \\
B^{(2,1)} & B^{(2,2)}
\end{pmatrix} = I - ADF(\bu),
\end{align*}
and split the norm between a finite number of columns and an estimate for the tail:
\begin{align*}
Z_1 &= \max\left( Z_1^{finite},\ Z_1^{tail} \right).
\end{align*}
$Z_1^{finite}$, which corresponds to
\begin{align*}
 \max_{0\leq n \leq 2N-2} \frac{1}{\xi_n(\nu)} \max\left(\Vert B^{(1,1)}_{(\cdot,n)} \Vert_\nu + \Vert B^{(2,1)}_{(\cdot,n)} \Vert_\nu,\ \Vert B^{(1,2)}_{(\cdot,n)} \Vert_\nu + \Vert B^{(2,2)}_{(\cdot,n)} \Vert_\nu\right),
\end{align*} 
can again be computed explicitly using interval arithmetic. A computation similar to the one performed in Section~\ref{sec:Z1_PM} shows that the second part, namely
\begin{align*}
\sup_{n\geq 2N-1} \frac{1}{\xi_n(\nu)} \max\left(\Vert B^{(1,1)}_{(\cdot,n)} \Vert_\nu + \Vert B^{(2,1)}_{(\cdot,n)} \Vert_\nu,\ \Vert B^{(1,2)}_{(\cdot,n)} \Vert_\nu + \Vert B^{(2,2)}_{(\cdot,n)} \Vert_\nu\right),
\end{align*}
can be bounded by
\begin{align}
\label{eq:Z1tail_SKT}
Z_1^{tail} =  \left\vert \left\Vert \begin{pmatrix} \1 & \0 \\ \0 & \1 
\end{pmatrix} - \bw \ast \Phi'(\bu) \right\Vert_\nu + \frac{1}{(N\pi)^2} \left\Vert \bw \right\Vert_\nu \left\Vert DR(\bu)\right\Vert_\nu \right\vert_1.
\end{align}

\begin{remark}
\label{rem:1_N2_SKT}
In the previous work~\cite{BreCas18} on computer-assisted proofs for the (triangular) SKT system, the estimate equivalent to $Z_1^{tail}$ scaled only as $\frac{1}{N}$, rather than as $\frac{1}{N^2}$ here (in practice the first term in~\eqref{eq:Z1tail_SKT} is small enough that only the second part really matters). Therefore, on top of not being restricted to the triangular case, our approach has the advantage of recovering the $\frac{1}{N^2}$ scaling expected for a second order equation with no first order term. This is comes from the fact that, contrarily to what was done in~\cite{BreCas18}, we do not have to resort to automatic differentiation techniques. This point will be discussed in more details in Section~\ref{sec:NP}.
\end{remark}

\subsubsection{The bound $Z_2$}

Again, we proceed as in Section~\ref{sec:Z2_PM} and estimate, for $\u$ and $\v$ in $\XX_\nu$,
\begin{align*}
\left\Vert AD^2F(\bu)(\u,\v) \right\Vert_{\XX_\nu} &\leq  \left(\Vert A^{(1,1)}\Delta \Vert_{\nu} + \Vert A^{(2,1)}\Delta \Vert_{\nu} \right) \left( 2d_{11}\Vert\u^{(1)}\Vert_\nu \Vert\v^{(1)}\Vert_\nu + d_{12}(\Vert\u^{(1)}\Vert_\nu \Vert\v^{(2)}\Vert_\nu + \Vert\v^{(1)}\Vert_\nu \Vert\u^{(2)}\Vert_\nu )\right)\\
& \quad + \left(\Vert A^{(1,1)} \Vert_{\nu} + \Vert A^{(2,1)} \Vert_{\nu} \right) \left( 2a_1\Vert\u^{(1)}\Vert_\nu \Vert\v^{(1)}\Vert_\nu + b_1(\Vert\u^{(1)}\Vert_\nu \Vert\v^{(2)}\Vert_\nu + \Vert\v^{(1)}\Vert_\nu \Vert\u^{(2)}\Vert_\nu )\right)\\
& \quad +  \left(\Vert A^{(1,2)}\Delta \Vert_{\nu} + \Vert A^{(2,2)}\Delta \Vert_{\nu} \right) \left( 2d_{22}\Vert\u^{(2)}\Vert_\nu \Vert\v^{(2)}\Vert_\nu + d_{21}(\Vert\u^{(1)}\Vert_\nu \Vert\v^{(2)}\Vert_\nu + \Vert\v^{(1)}\Vert_\nu \Vert\u^{(2)}\Vert_\nu )\right)\\
& \quad + \left(\Vert A^{(1,2)} \Vert_{\nu} + \Vert A^{(2,2)} \Vert_{\nu} \right) \left( 2a_2\Vert\u^{(2)}\Vert_\nu \Vert\v^{(2)}\Vert_\nu + b_2(\Vert\u^{(1)}\Vert_\nu \Vert\v^{(2)}\Vert_\nu + \Vert\v^{(1)}\Vert_\nu \Vert\u^{(2)}\Vert_\nu )\right)\\
&\leq \max\left[ 2d_{11} \left(\Vert A^{(1,1)}\Delta \Vert_{\nu} + \Vert A^{(2,1)}\Delta \Vert_{\nu}\right) + 2a_1 \left(\Vert A^{(1,1)} \Vert_{\nu} + \Vert A^{(2,1)} \Vert_{\nu}\right) , \right. \\
&\qquad\qquad  2d_{22} \left(\Vert A^{(1,2)}\Delta \Vert_{\nu} + \Vert A^{(2,2)}\Delta \Vert_{\nu}\right) + 2a_2 \left(\Vert A^{(1,2)} \Vert_{\nu} + \Vert A^{(2,2)} \Vert_{\nu}\right) ,  \\
&\qquad\qquad  d_{12} \left(\Vert A^{(1,1)}\Delta \Vert_{\nu} + \Vert A^{(2,1)}\Delta \Vert_{\nu}\right) + b_1 \left(\Vert A^{(1,1)} \Vert_{\nu} + \Vert A^{(2,1)}\Vert_{\nu}\right) \\
&\qquad\qquad + \left. d_{21} \left(\Vert A^{(1,2)}\Delta \Vert_{\nu} + \Vert A^{(2,2)}\Delta \Vert_{\nu}\right) + b_2 \left(\Vert A^{(1,2)} \Vert_{\nu} + \Vert A^{(2,2)} \Vert_{\nu}\right) \right] \left\Vert \u\right\Vert_{\XX_\nu} \left\Vert \v\right\Vert_{\XX_\nu},
\end{align*}
which yields
\begin{align*}
Z_2 &= \max\left[ 2d_{11} \left(\Vert A^{(1,1)}\Delta \Vert_{\nu} + \Vert A^{(2,1)}\Delta \Vert_{\nu}\right) + 2a_1 \left(\Vert A^{(1,1)} \Vert_{\nu} + \Vert A^{(2,1)} \Vert_{\nu}\right) , \right. \\
&\qquad\qquad  2d_{22} \left(\Vert A^{(1,2)}\Delta \Vert_{\nu} + \Vert A^{(2,2)}\Delta \Vert_{\nu}\right) + 2a_2 \left(\Vert A^{(1,2)} \Vert_{\nu} + \Vert A^{(2,2)} \Vert_{\nu}\right) ,  \\
&\qquad\qquad  d_{12} \left(\Vert A^{(1,1)}\Delta \Vert_{\nu} + \Vert A^{(2,1)}\Delta \Vert_{\nu}\right) + b_1 \left(\Vert A^{(1,1)} \Vert_{\nu} + \Vert A^{(2,1)}\Vert_{\nu}\right) \\
&\qquad\qquad + \left. d_{21} \left(\Vert A^{(1,2)}\Delta \Vert_{\nu} + \Vert A^{(2,2)}\Delta \Vert_{\nu}\right) + b_2 \left(\Vert A^{(1,2)} \Vert_{\nu} + \Vert A^{(2,2)} \Vert_{\nu}\right) \right].
\end{align*}

\subsection{Examples and results}
\label{sec:res_SKT}

We now consider several specific parameter values, and use the presented methodology to obtain qualitative existence theorems about the steady states of the SKT system~\eqref{eq:SKT} for those parameter values, with $\Omega = (0,1)$. The different parameter sets that we study here are given in Table~\ref{tab:paraSKT}, and their respective relevance is explained in the following subsections.

\begin{table}[!h]
\centering
\begin{tabular}{>{\columncolor[gray]{0.9}}c>{\columncolor[gray]{0.9}}c>{\columncolor[gray]{0.9}}c>{\columncolor[gray]{0.9}}c>{\columncolor[gray]{0.9}}c>{\columncolor[gray]{0.9}}cccccc c>{\columncolor[gray]{0.7}}c}
$d_1$&$d_2$&$d_{12}$&$d_{21}$&$d_{11}$&$d_{22}$&$r_1$&$r_2$&$a_1$&$a_2$&$b_1$&$b_2$&\\
\midrule[2pt]
0.005 & 0.005 & 3 & 0 & 0 & 0 & 5 & 2 & 3 & 3 & 1 & 1 & Section~\ref{sec:SKT1} \\
\midrule
0.005 & 0.005 & 100 & 100 & 0 & 0 & 15/2 & 16/7 & 4 & 2 & 6 & 1 & Section~\ref{sec:SKT2} \\
\midrule
0.05 & 0.05 & 3 & 0 & 0 & 0 & 15 & 5 & 1 & 1 & 0.5 & 3 & Section~\ref{sec:SKT3} \\
\midrule
-0.007 & -0.007 & 3 & 0.002 & 0.05 & 0.05 & 5 & 2 & 3 & 3 & 1 & 1 & Section~\ref{sec:SKT4} \\
\midrule
\bottomrule[2pt]
\end{tabular}
\caption{The parameter sets of the SKT system~\eqref{eq:SKT} considered in this work.}
\label{tab:paraSKT}
\end{table}

In each case, we only show one or a couple of steady states, but we emphasize that there might exist many more steady states for the exact same parameter values. We refer to~\cite{IidMimNim06,IzuMim08,BreKueSor21} for a broader picture and many bifurcation diagrams.

Finally, before getting to each case let us also mention that, even if we only study the existence and the precise description of (some of) the steady states here, their stability is also a very important question. While adapting the computer-assisted techniques presented in this work to prove that a steady steady is unstable is relatively straightforward (see e.g.~\cite{BreCas18} for the triangular case), a computer-assisted approach that could be used to prove (linear) stability for an arbitrary solution of the SKT system is still lacking.   

\subsubsection{Comparison with the previous setup using automatic differentiation~\cite{BreCas18}}
\label{sec:SKT1}

We first focus on the first parameter set (i.e. the first row) in Table~\ref{tab:paraSKT}. For these parameter values, a complex bifurcation diagram of steady states was first computed in~\cite{IidMimNim06}, and then mostly validated in~\cite{BreCas18} using an approach relying on the fact that $d_{11}=d_{22}=d_{21}=0$. We therefore use this first case as a benchmark, and compare the results of our method with the ones of~\cite{BreCas18}.  

\begin{figure}[h!]
\centering
\includegraphics[scale=0.5]{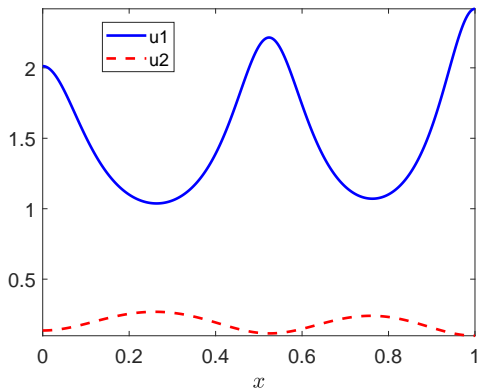}
\caption{The approximate steady state $\bar u$ of the SKT systsem~\eqref{eq:SKT} on $\Omega=(0,1)$ and with parameters as in the first row of Table~\ref{tab:paraSKT}, which has been validated in Theorem~\ref{th:SKT1}.}
\label{fig:SKT1}
\end{figure}

\begin{theorem}
\label{th:SKT1}
Let $\bar u$ be the function represented on Figure~\ref{fig:SKT1}, and whose precise description in terms of Fourier coefficients can be found at~\cite{Bre21}. There exists a smooth steady state $u$ of the SKT system~\eqref{eq:SKT}, with parameters as in the first row of Table~\ref{tab:paraSKT} and $\Omega=(0,1)$, such that $\left\Vert u-\bar u\right\Vert_{C^0} \leq 4\times 10^{-7}$.
\end{theorem}
\begin{proof}
We take $N=50$ and $\nu=1.01$, compute the estimates $Y$, $Z_1$ and $Z_2$ obtained in Section~\ref{sec:bounds_SKT}, the approximate values for the bounds being
\begin{align*}
Y = 2.7\times 10^{-7},\quad Z_1 = 0.25,\quad Z_2 = 14000,
\end{align*} 
and apply Theorem~\ref{th:NK_SKT}. 
The computational parts of the proof can be reproduced by running the Matlab code \texttt{script\_SKT.m} from~\cite{Bre21}, together with Intlab~\cite{Rum99}.
\end{proof}

Remark~\ref{rem:space} also applies here (and for the remainder of Section~\ref{sec:res_SKT}). The main point that we want to emphasize with this theorem is not the result itself, which was already obtained in~\cite{BreCas18}, but the fact that our new approach is significantly more efficient. Indeed, with the same parameters ($N=50$ and $\nu=1.01$) the estimate corresponding to $Z_1^{tail}$ in~\cite{BreCas18} would be approximately equal to $1.46\ldots$ and therefore assumption~\eqref{eq:Z1l1} would not hold. Comparatively, with our approach we get $Z_1^{tail} \approx 0.25$. This means the validation is less costly with our approach, because it works with $N$ smaller, and this comes from the fact that we avoid the usage of automatic differentiation techniques in (see Remark~\ref{rem:1_N2_SKT}). Moreover, with $N=50$, $Z_1$ is still significantly away from $1$ in our approach, and the main reason why we need close to $N=50$ modes to validate this solution is only because otherwise $Y$ is not small enough. This is exactly what one would hope for, i.e. that the validation is successful as soon as we have enough modes for the approximate solution to be sufficiently accurate.

\subsubsection{A first non-triangular case}
\label{sec:SKT2}

We now focus on the second parameter set in Table~\ref{tab:paraSKT}. Notice that we are no longer in the so-called triangular case (both $d_{12}$ and $d_{21}$ are non zero), therefore the previous approach from~\cite{BreCas18} can no longer be used. Moreover, the solutions which are proved to exist for these parameters are of particular interest, because they answer positively a question which was open until recently, about the existence of non homogeneous steady states in the weak competition case ($b_1b_2<a_1a_2$) when both cross diffusion coefficients are large, see~\cite{BreKueSor21} and the reference therein for more details. In~\cite{BreKueSor21} this question was answered positively by a bifurcation analysis which provides information locally (close to the bifurcation point), and solutions far away from the bifurcation were computed numerically. It is such solutions far away from the bifurcation whose existence we are able to prove here with computer-assistance.

\begin{figure}[h!]
\centering
\begin{subfigure}[t]{0.49\textwidth}
\centering
\includegraphics[scale=0.49]{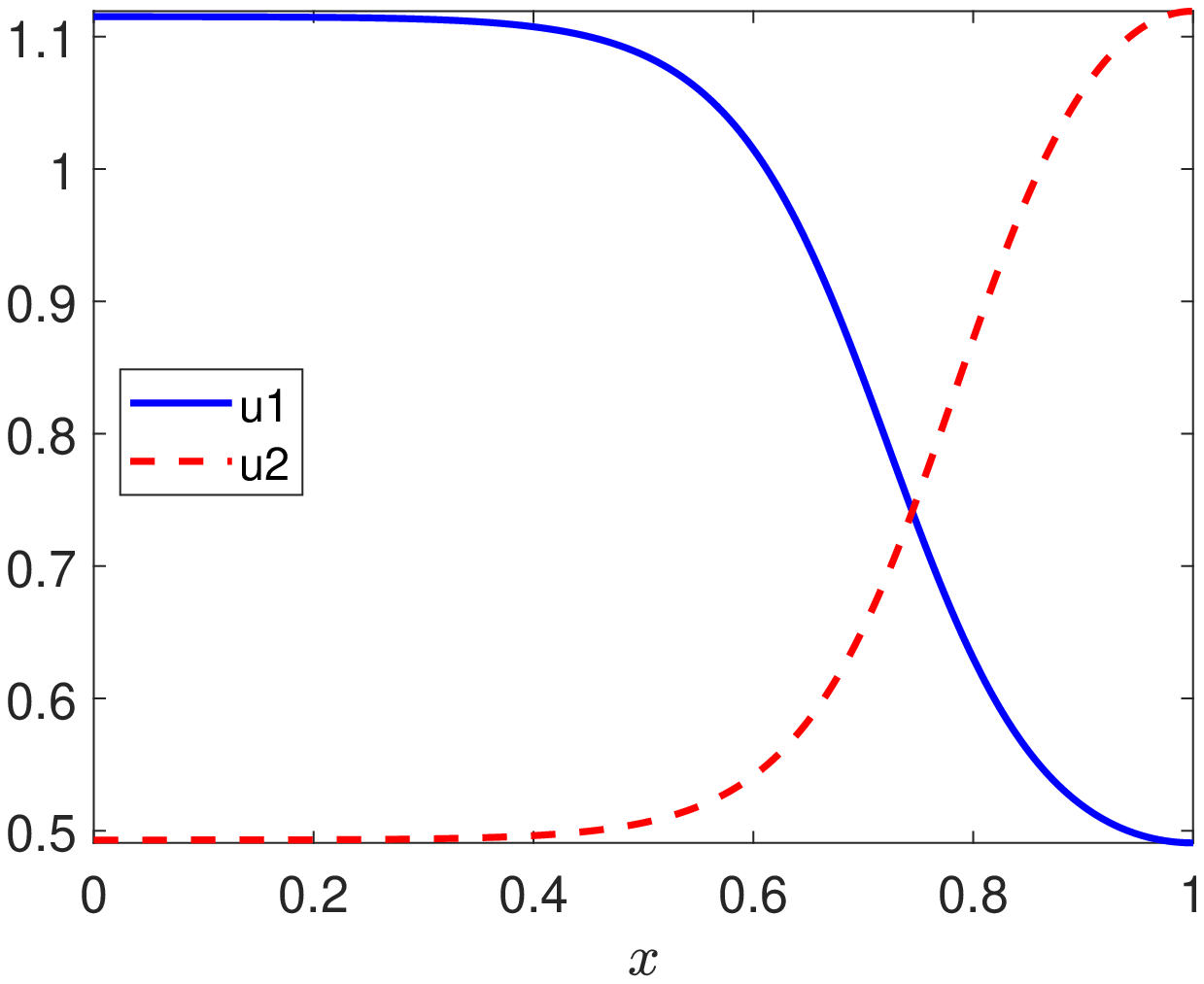}
\end{subfigure}
\begin{subfigure}[t]{0.49\textwidth}
\centering
\includegraphics[scale=0.49]{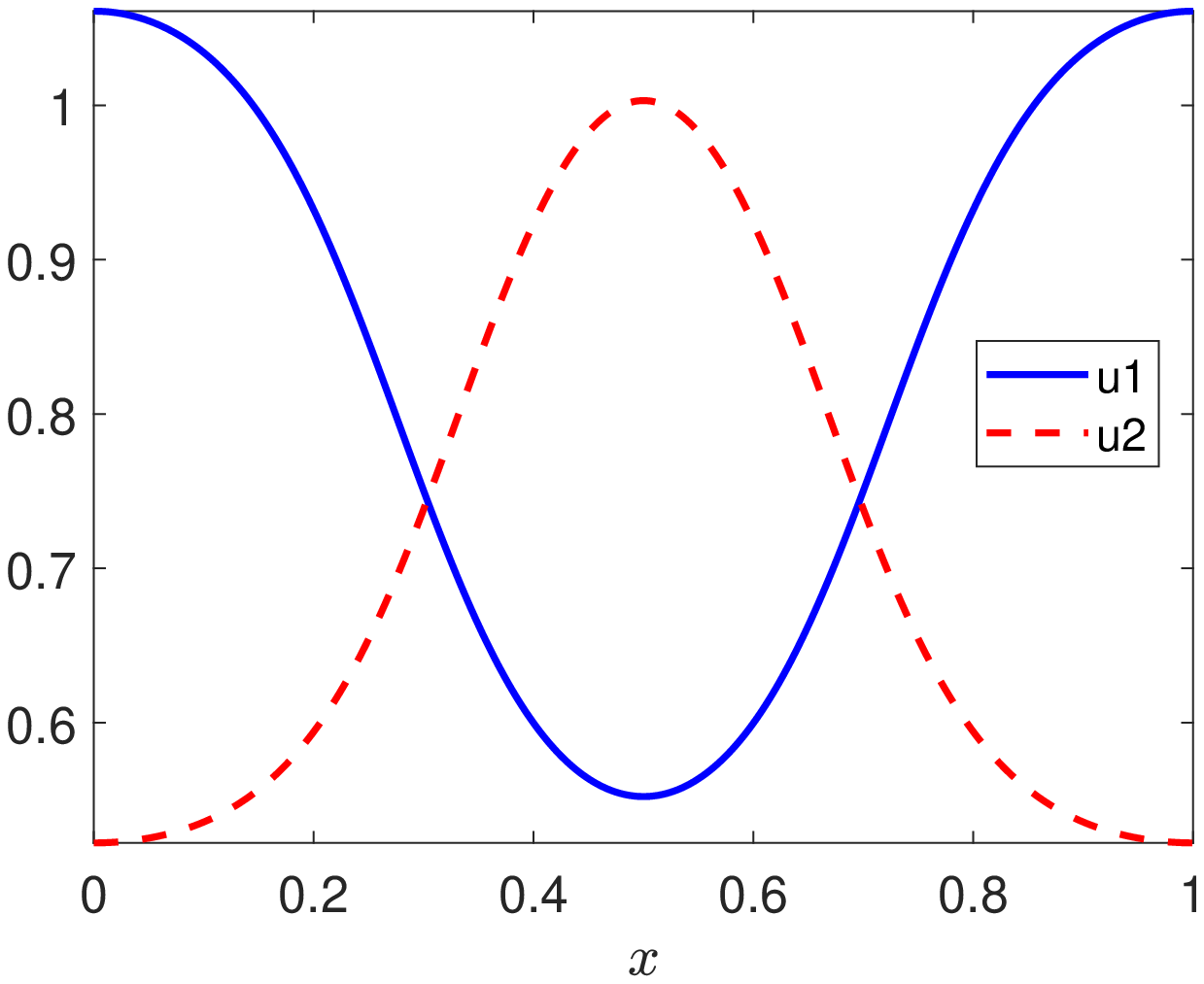}
\end{subfigure}
\caption{The two approximate steady states of the SKT systsem~\eqref{eq:SKT} on $\Omega=(0,1)$ and with parameters as in the second row of Table~\ref{tab:paraSKT}, which have been validated in Theorem~\ref{th:SKT2}.}
\label{fig:SKT2}
\end{figure}

\begin{theorem}
\label{th:SKT2}
For each of the functions represented on Figure~\ref{fig:SKT2}, and whose precise description in terms of Fourier coefficients can be found at~\cite{Bre21}, there exists a smooth steady state of the SKT system~\eqref{eq:SKT}, with parameters as in the second row of Table~\ref{tab:paraSKT} and $\Omega=(0,1)$, at a distance of at most $8\times 10^{-12}$ in $C^0$ norm.
\end{theorem}
\begin{proof}
We take $N=50$ and $\nu=1.01$, compute the estimates $Y$, $Z_1$ and $Z_2$ obtained in Section~\ref{sec:bounds_SKT} for each approximate solution, and apply Theorem~\ref{th:NK_SKT}. 
The computational parts of the proof can be reproduced by running the Matlab code \texttt{script\_SKT.m} from~\cite{Bre21}, together with Intlab~\cite{Rum99}, which provides more details such as the precise value of each bound for each solution.
\end{proof}

\subsubsection{A new parameter regime}
\label{sec:SKT3}

We now focus on the third parameter set in Table~\ref{tab:paraSKT}. In order to understand it's interest, it is helpful to briefly consider the homogeneous case, i.e.~\eqref{eq:SKT} with only the reaction terms, for which one can easily identify three different parameter regimes with different dynamics, see e.g.~\cite{LouNi96}.
\begin{itemize}
\item \textbf{Case 1 (weak competition):} 
\begin{align*}
\frac{b_1}{a_2} < \frac{r_1}{r_2} < \frac{a_1}{b_2}.
\end{align*}
Here, any solution with positive initial data converges to a co-existence steady state 
\begin{align}
\label{eq:SKTeq}
\left( \frac{r_1a_2-r_2b_1}{a_2a_1-b_2b_1},\frac{r_2a_1-r_1b_2}{a_2a_1-b_2b_1}\right).
\end{align}
\item \textbf{Case 2 (strong competition):} 
\begin{align*}
\frac{a_1}{b_2} < \frac{r_1}{r_2} < \frac{b_1}{a_2}.
\end{align*}
Here, the above co-existence steady state~\eqref{eq:SKTeq} is still positive, but it is unstable, and generically a solution with positive initial data converges to one of the two extinction states
\begin{align}
\label{eq:SKText}
\left( \frac{r_1}{a_1}, 0\right) \qquad \text{or}\qquad \left(0 , \frac{r_2}{a_2}\right),
\end{align}
depending on the initial condition.
\item \textbf{Case 3:}
everything else (except the case $b_1/a_2 = r_1/r_2 = a_1/b_2$ which is degenerate). Here, any solution with positive initial data also converges to one of the two extinction states~\eqref{eq:SKText}, but all solutions converge to the same one (the one being selected depends on the parameter values, not on the initial data).
\end{itemize}
The aim of introducing a spatial dimension to model, and more precisely nonlinear diffusion as in~\eqref{eq:SKT}, is to obtain more interesting and realistic solutions, and in particular non-homogeneous positive steady states. However, it is striking that in the literature on the SKT system, the existence of non-homogeneous steady states is usually studied in cases 1 and 2 only. This might be related to the fact that the analysis in case 3 can be more difficult, because there is no positive steady state that can be used as a kind of starting point: the co-existence steady states~\eqref{eq:SKTeq} still exists (except if $a_1a_2=b_1b_2$), but it is no longer positive in case 3.

However, recent numerical experiments suggested that positive non-homogeneous do also exist in case 3~\cite{BreKueSor21}. The third parameter set in Table~\ref{tab:paraSKT} corresponds to such a case, and we now rigorously establish the existence of (some of) the solutions that were found numerically in~\cite{BreKueSor21}.

\begin{figure}[h!]
\centering
\begin{subfigure}[t]{0.49\textwidth}
\centering
\includegraphics[scale=0.49]{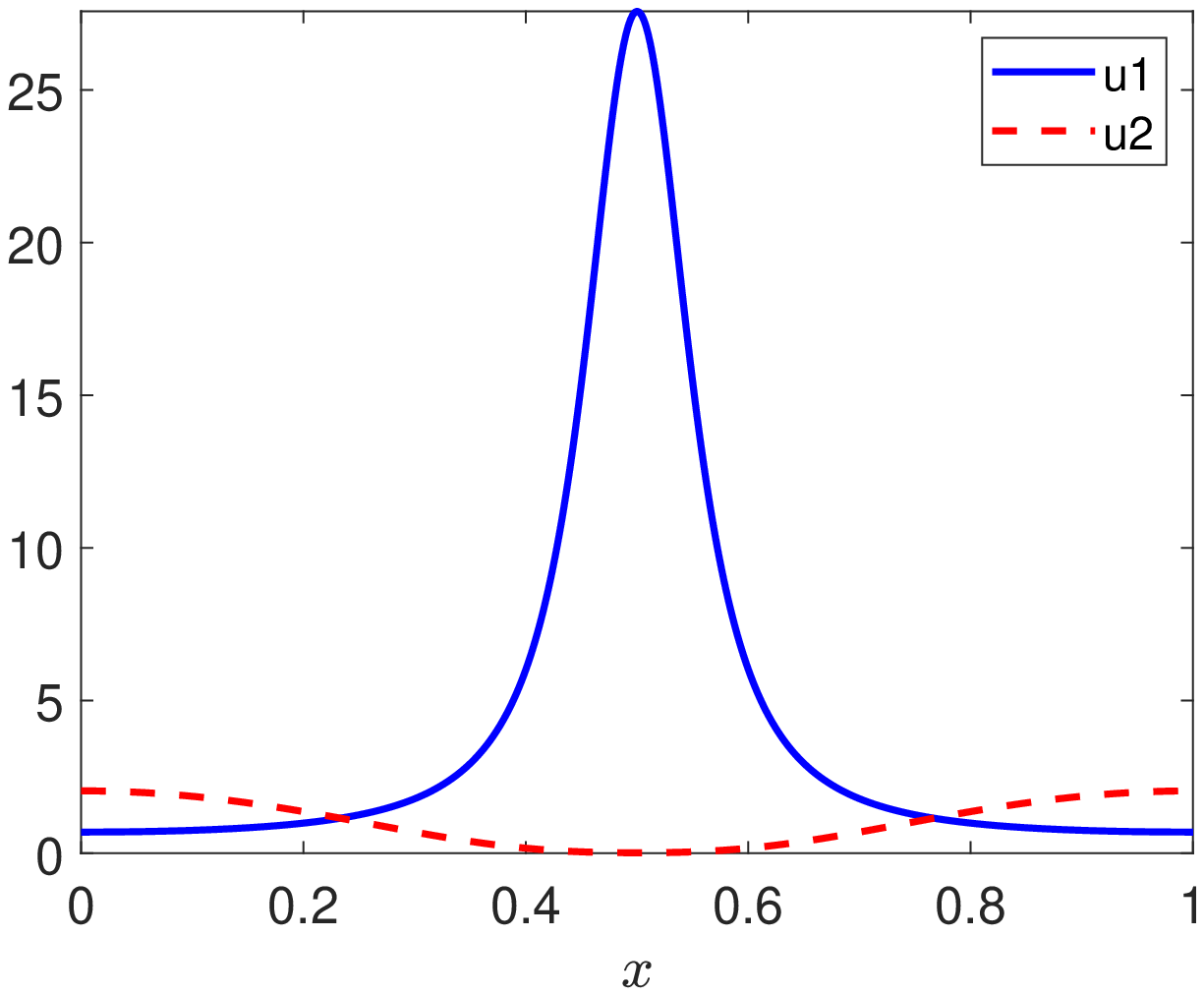}
\end{subfigure}
\begin{subfigure}[t]{0.49\textwidth}
\centering
\includegraphics[scale=0.49]{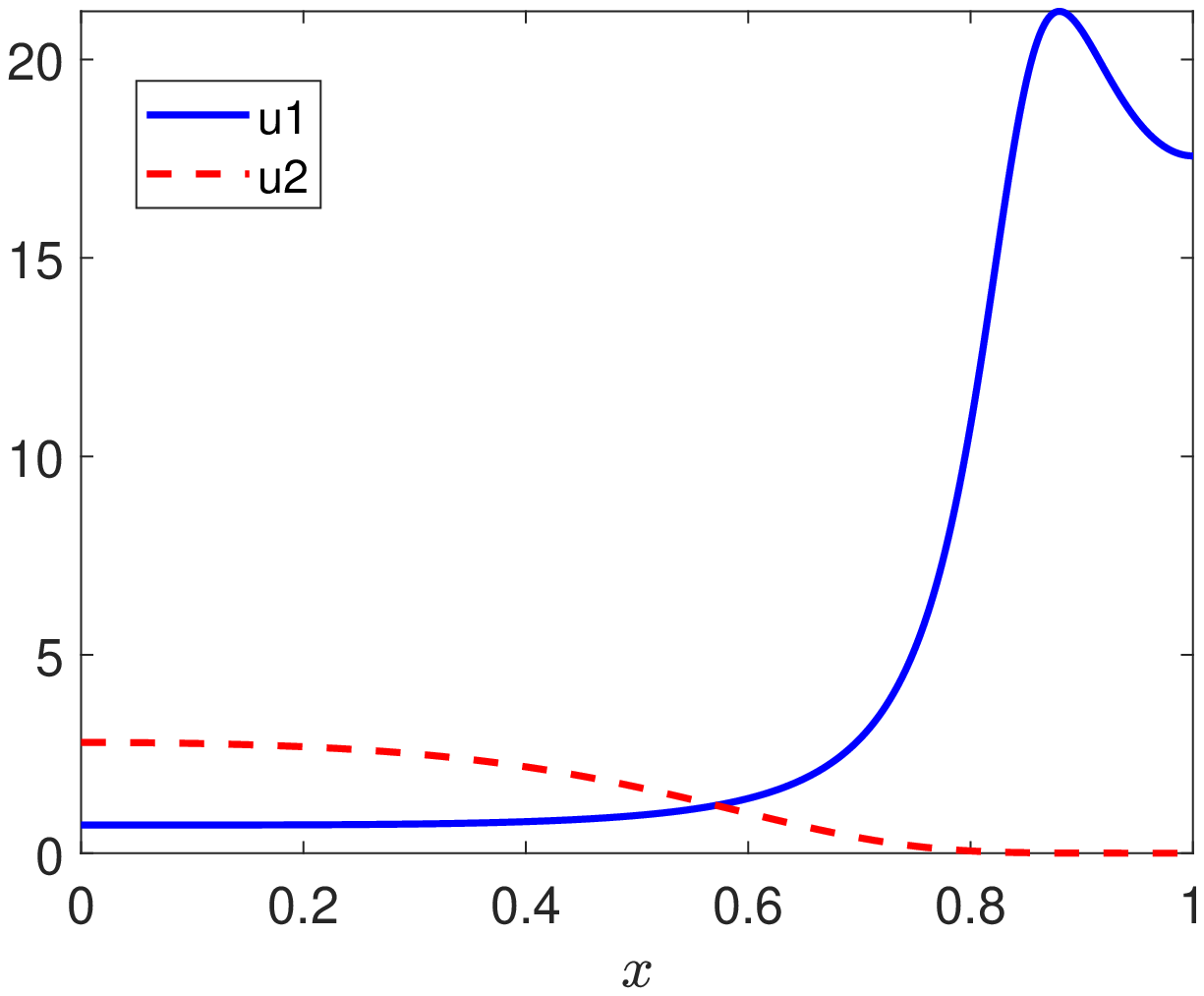}
\end{subfigure}
\caption{The two approximate steady states of the SKT systsem~\eqref{eq:SKT} on $\Omega=(0,1)$ and with parameters as in the third row of Table~\ref{tab:paraSKT}, which have been validated in Theorem~\ref{th:SKT3}.}
\label{fig:SKT3}
\end{figure}

\begin{theorem}
\label{th:SKT3}
For each of the functions represented on Figure~\ref{fig:SKT3}, and whose precise description in terms of Fourier coefficients can be found at~\cite{Bre21}, there exists a smooth steady state of the SKT system~\eqref{eq:SKT}, with parameters as in the third row of Table~\ref{tab:paraSKT} and $\Omega=(0,1)$, at a distance of at most $9\times 10^{-11}$ in $C^0$ norm.
\end{theorem}
\begin{proof}
We take $N=500$ and $\nu=1.005$,  compute the estimates $Y$, $Z_1$ and $Z_2$ obtained in Section~\ref{sec:bounds_SKT} for each approximate solution, and apply Theorem~\ref{th:NK_SKT}.
The computational parts of the proof, can be reproduced by running the Matlab code \texttt{script\_SKT.m} from~\cite{Bre21}, together with Intlab~\cite{Rum99}, which provides more details such as the precise value of each bound for each solution.
\end{proof}

\subsubsection{Full case}
\label{sec:SKT4}

We finally focus on the fourth and last parameter set in Table~\ref{tab:paraSKT}, which allows us to showcase that our approach is still successful when all coefficients related to the diffusion term are non zero. As a bonus, we show that positive self-diffusion coefficients $d_{11}$ and $d_{22}$ allow to take negative values of the linear diffusion coefficients $d_1$ and $d_2$ and still retain positive and smooth solutions. This was also already noticed numerically in~\cite{BreKueSor21}.

\begin{figure}[h!]
\centering
\begin{subfigure}[t]{0.49\textwidth}
\centering
\includegraphics[scale=0.49]{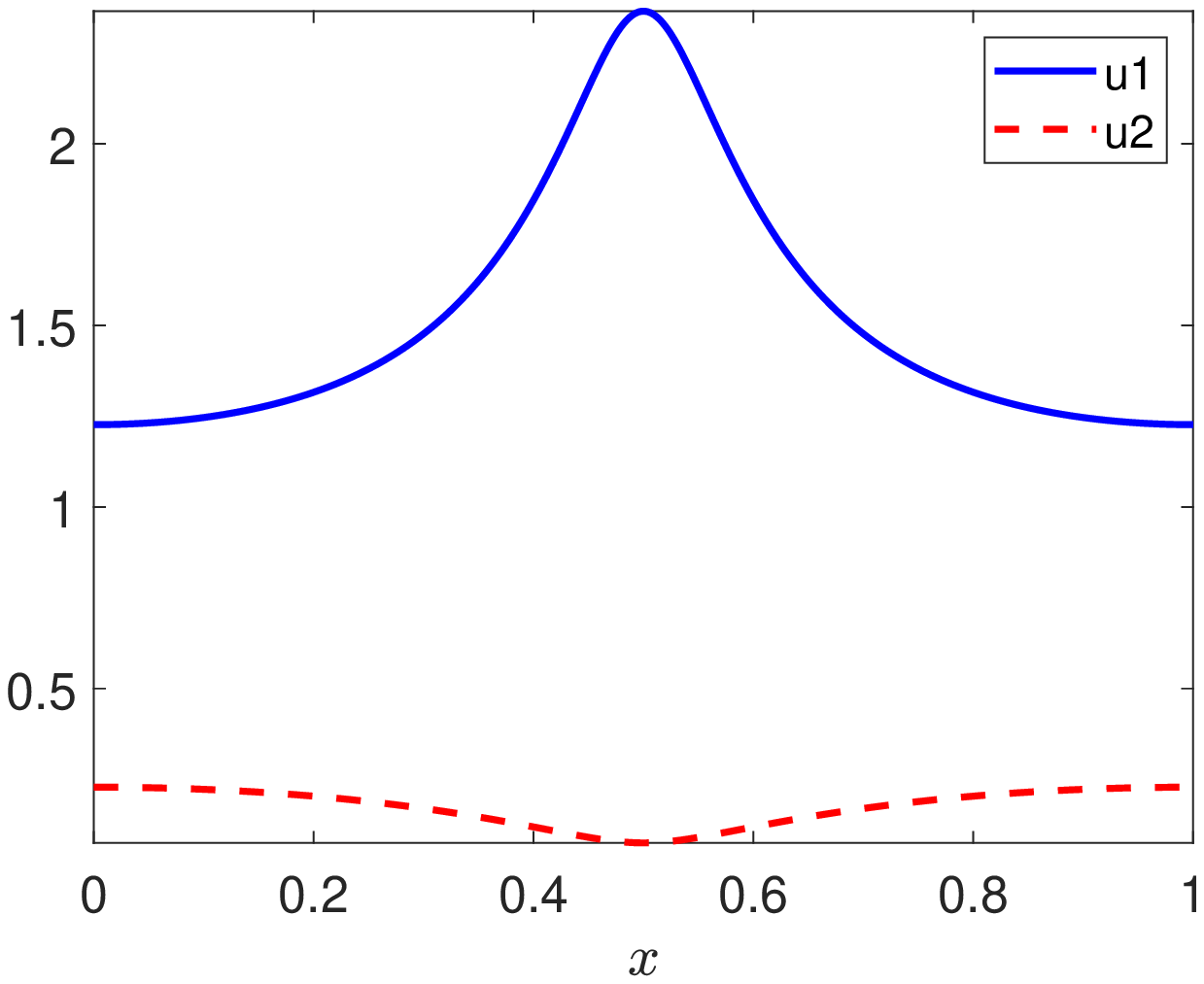}
\end{subfigure}
\begin{subfigure}[t]{0.49\textwidth}
\centering
\includegraphics[scale=0.49]{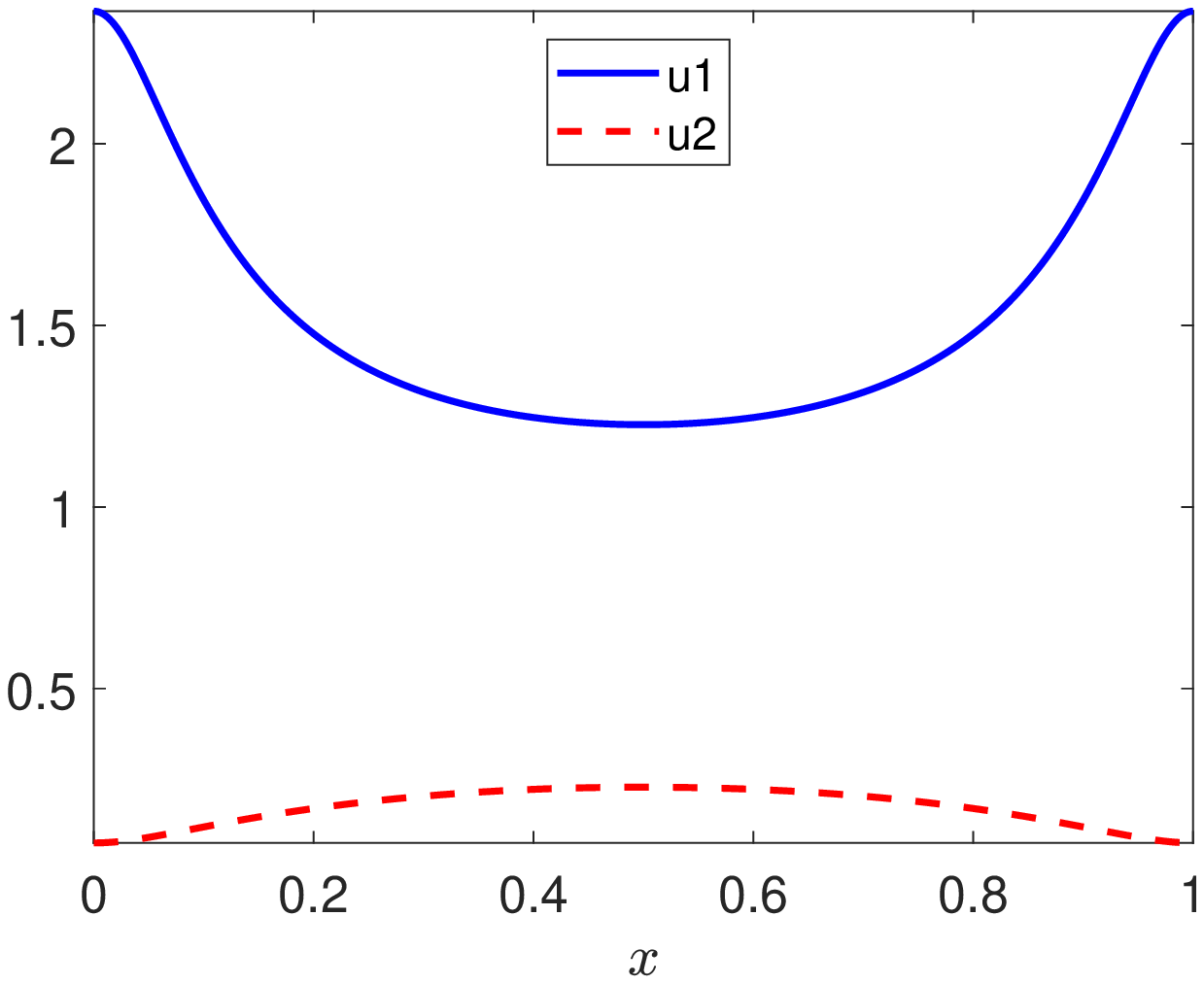}
\end{subfigure}
\begin{subfigure}[t]{0.49\textwidth}
\centering
\includegraphics[scale=0.49]{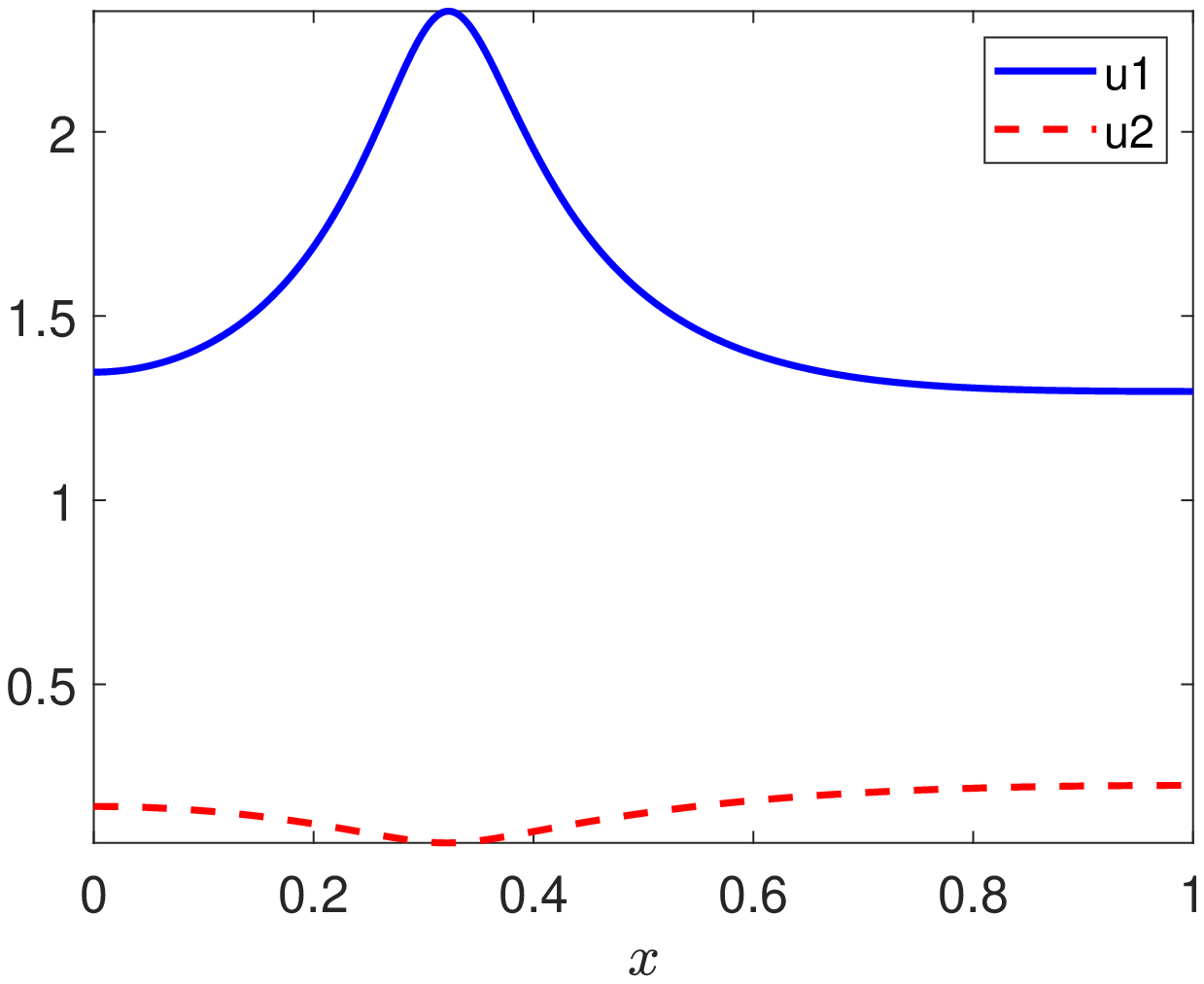}
\end{subfigure}
\begin{subfigure}[t]{0.49\textwidth}
\centering
\includegraphics[scale=0.49]{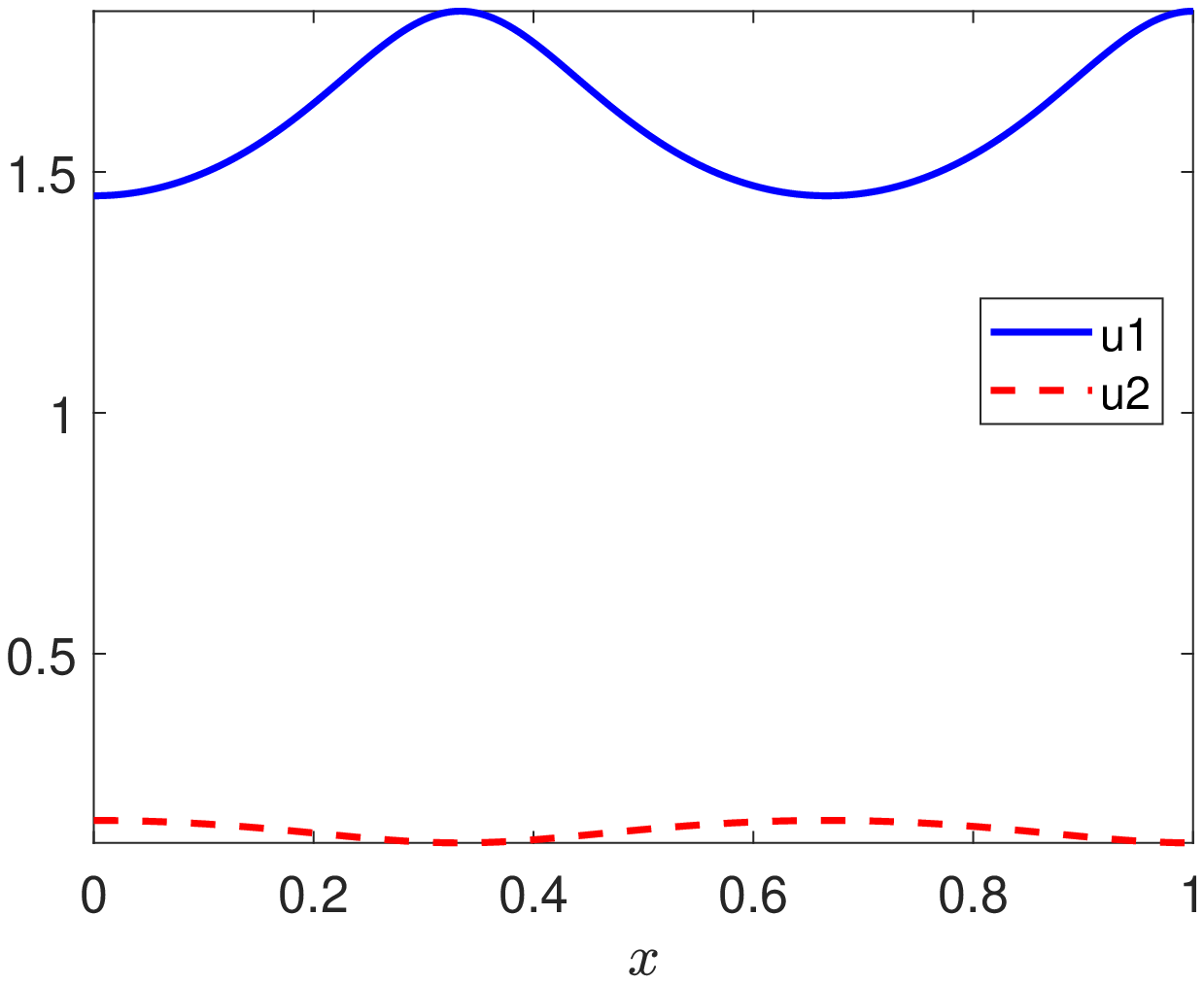}
\end{subfigure}
\caption{The four approximate steady states of the SKT systsem~\eqref{eq:SKT} on $\Omega=(0,1)$ and with parameters as in the fourth row of Table~\ref{tab:paraSKT}, which have been validated in Theorem~\ref{th:SKT4}.}
\label{fig:SKT4}
\end{figure}

\begin{theorem}
\label{th:SKT4}
For each of the functions represented on Figure~\ref{fig:SKT4}, and whose precise description in terms of Fourier coefficients can be found at~\cite{Bre21}, there exists a smooth steady state of the SKT system~\eqref{eq:SKT}, with parameters as in the fourth row of Table~\ref{tab:paraSKT} and $\Omega=(0,1)$, at a distance of at most $1\times 10^{-9}$ in $C^0$ norm.
\end{theorem}
\begin{proof}
We take $N=100$ and $\nu=1.01$,  compute the estimates $Y$, $Z_1$ and $Z_2$ obtained in Section~\ref{sec:bounds_SKT} for each approximate solution, and apply Theorem~\ref{th:NK_SKT}.
The computational parts of the proof, can be reproduced by running the Matlab code \texttt{script\_SKT.m} from~\cite{Bre21}, together with Intlab~\cite{Rum99}, which provides more details such as the precise value of each bound for each solution.
\end{proof}

\section{Beyond polynomial diffusion terms}
\label{sec:NP}

In this section, we go back to a scalar problem on a one-dimensional domain of the form
\begin{align}
\label{eq:PMbis}
\left\{
\begin{aligned}
&\Delta \Phi(u) + R(u) = 0 \qquad & \text{on }(0,1), \\
&\frac{\partial u}{\partial n } = 0 \qquad & \text{on } \{0,1\},
\end{aligned}
\right.
\end{align}
with $R(u) = \alpha u -\beta u^2 + g$, but this time with a non-polynomial term $\Phi$ of the form $\Phi(u)=\frac{u}{\gamma + u}$. Again, rather than aiming for the most general result possible, we chose to consider an explicit scalar example to simplify the presentation, but the ideas presented in this section can also be applied with different non-polynomial terms, systems, and higher space dimension.

If the solution $u$ under consideration satisfies $\Vert \u\Vert_\nu < \gamma$, then $\Phi(\u) = \sum_{k\geq 0} (-1)^k \frac{\u^{k+1}}{\gamma^{k+1}}$ is still well defined in $\ell^1_\nu$, and the procedure described in Section~\ref{sec:PM} could be adapted, at the expense of more pen and paper estimates. However, as we are going to see, there are solutions for which $\Vert \u\Vert_\nu > \gamma$.

Even in the case where $\Vert \u\Vert_\nu < \gamma$, it might be more convenient to try to get back to a polynomial system, by introducing new variables. Following the ideas of~\cite{LesMirRan16}, one could consider the additional variable $v= \frac{u}{\gamma + u}$, rewrite~\eqref{eq:PMbis} as a first order system, and append to it the differential equation satisfied by $v$. While this approach can be appealing, because the resulting system is again polynomial and therefore fits immediately into the existing framework, it also has a couple of shortcomings. First, having to go back to a system of first order equations is a downside, because all the $\frac{1}{N^2}$ terms in the estimates (like in~\eqref{eq:Z1tail_PM} and~\eqref{eq:Z1tail_SKT}) are now replaced by $\frac{1}{N}$ terms, and therefore the dimension $N$ of the projection has to be taken much larger for these terms to be small enough. Second, it does not seem straightforward to generalize these ideas as soon as the space dimension or the number of unknown is greater than one.

In order to overcome these limitations, we propose a new approach, which is based on the following observation: Introducing new variables is a good idea, but in the context of computer-assisted proofs there is nothing wrong with differential-algebraic equations (especially when the algebraic equations can be rewritten as polynomial ones), and in particular there is no need to differentiate the algebraic equations. That is, we do also introduce the additional variable $v= \frac{u}{\gamma + u}$, but we then directly work on the system
\begin{align}
\label{eq:NP}
\left\{
\begin{aligned}
&\Delta v + R(u) = 0 \qquad & \text{on }(0,1),  \\
&u - (\gamma + u)v = 0 \qquad & \text{on }(0,1), \\
&\frac{\partial u}{\partial n } = 0 = \frac{\partial u}{\partial n } \qquad & \text{on }\{0,1\}.
\end{aligned}
\right.
\end{align}
Indeed, we show in this section that this system of differential-algebraic equations is amenable to a posteriori validation techniques very similar to the ones presented up to now in this paper.

\subsection{$F=0$ into fixed-point problem}
\label{sec:fp_NP}

Relabeling the unknowns $(u,v)$ in~\eqref{eq:NP} as $(u^{(1)},u^{(2)})$, we can directly use all the notations, spaces and norms introduced in Section~\ref{sec:SKT}. That is, we look for a zero $\u=(\u^{(1)},\u^{(2)})$ in $\XX_\nu$ of
\begin{align*}
F=\begin{pmatrix}
F^{(1)} \\
F^{(2)}
\end{pmatrix} 
\end{align*}
defined as
\begin{align*}
\left\{
\begin{aligned}
F^{(1)}_n(\u) &= \u^{(1)}_n -\gamma \u^{(2)}_n - (\u^{(1)}\ast\u^{(2)})_n \\
F^{(2)}_n(\u) &= -(n\pi)^2\u^{(2)}_n + R_n(\u^{(1)})
\end{aligned}
\right.
\qquad \forall~n\geq 0,
\end{align*}
or, in a more condensed form,
\begin{align*}
F(\u) = \begin{pmatrix}
\u^{(1)} -\gamma\u^{(2)} - \u^{(1)}\ast\u^{(2)}\\
\Delta \u^{(2)} + R(\u^{(1)})
\end{pmatrix},
\end{align*}
with $R(u) = \alpha u -\beta u^2 + g$.

We now assume that we have computed an approximate zero $\bu$ in $\Pi_N\XX_\nu$ that we want to validate a posteriori. We also assume that we have computed $\bw\in\Pi_N\ell^1_\nu$ such that $\bw \ast (\1-\bu^{(2)}) \approx \1$, and $\bs\in\Pi_N\ell^1_\nu$ such that $\bs \approx \bw\ast(\bm{\gamma}+\bu^{(1)})$.
Next, we consider 
\begin{align*}
\bar A = 
\renewcommand{\arraystretch}{2}
\left(
\begin{array}{c|c}
\bar A^{(1,1)} & \bar A^{(1,2)} \\ \hline
\bar A^{(2,1)} & \bar A^{(2,2)}
\end{array}
\right),
\end{align*}
a numerically computed approximation of $\Pi_N \left ( DF(\bu)^{-1}\right) \Pi_N $ (Remark~\ref{rem:barA} is also relevant here). The approximate inverse $A$ for this problem is then defined as
\begin{align}
\label{eq:A_NP}
\newcommand{\phs}{\phantom{\bar A}}
\newcommand{\phl}{\phantom{M(\sigma)\Delta^{-1}}}
A = 
\left(
\begin{array}{c|c}
\begin{array}{cccccc}
 & & & \multicolumn{1}{|c}{ } &  & \\ 
 & \bar A^{(1,1)} & & \multicolumn{1}{|c}{ } &  & \\ 
\phs & \phs & \phs & \multicolumn{1}{|c}{ } &  & \\ \cline{1-3}
 & & & & & \\
 & & & & M(\bw) & \\
 & & & & \phl &
\end{array}
&
\begin{array}{cccccc}
 & & & \multicolumn{1}{|c}{ } &  & \\ 
 & \bar A^{(1,2)} & & \multicolumn{1}{|c}{ } &  & \\ 
\phs & \phs & \phs & \multicolumn{1}{|c}{ } &  & \\ \cline{1-3}
 & & & & & \\
 & & & & M(\bs)\Delta^{-1} & \\
 & & & & &
\end{array} \\ \hline
\begin{array}{cccccc}
 & & & \multicolumn{1}{|c}{ } &  & \\ 
 & \bar A^{(2,1)} & & \multicolumn{1}{|c}{ } &  & \\ 
\phs & \phs & \phs & \multicolumn{1}{|c}{ } &  & \\ \cline{1-3}
 & & & & & \\
 & & & & 0 & \\
 & & & & \phl &
\end{array}
&
\begin{array}{cccccc}
 & & & \multicolumn{1}{|c}{ } &  & \\ 
 & \bar A^{(2,2)} & & \multicolumn{1}{|c}{ } &  & \\ 
\phs & \phs & \phs & \multicolumn{1}{|c}{ } &  & \\ \cline{1-3}
 & & & & & \\
 & & & & \Delta^{-1} & \\
 & & & & \phl &
\end{array}
\end{array}
\right),
\end{align}
or more compactly as
\begin{align*}
A= \bar A + \left(\tilde A - \Pi_N  \tilde A \Pi_N \right),
\end{align*}
where  $\tilde A$ is defined by
\begin{align*}
\tilde A \u := \begin{pmatrix}
\bw \ast \u^{(1)} + \bs \ast \Delta^{-1}\u^{(2)} \\
\Delta^{-1}\u^{(1)}
\end{pmatrix}.
\end{align*}

Now that $A$ is defined, the validation conditions are again the same as in Theorem~\ref{th:NK} and~\ref{th:NK_SKT}, the only difference being the extra condition ensuring the injectivity of $A$, which is differs because the structure of $A$ does.
\begin{theorem}
\label{th:NK_NP}
With the notations introduced in this section, assume there exist constants $Y$, $Z_1$ and $Z_2$ satisfying
\begin{align*}
\Vert AF(\bu) \Vert_{\XX_\nu} &\leq Y \\
\Vert I - ADF(\bu) \Vert_{\XX_\nu}  &\leq Z_1 \\
\Vert AD^2F(\u) \Vert_{\XX_\nu} &\leq Z_2 \qquad \forall~\u\in\XX_\nu, 
\end{align*}
and
\begin{align*}
Z_1 &< 1 \\
2YZ_2 &< (1-Z_1)^2.
\end{align*}
Then, for any $r$ satisfying 
\begin{align*}
\frac{1-Z_1 - \sqrt{(1-Z_1)^2-2YZ_2}}{Z_2} \leq r < \frac{1-Z_1}{Z_2},
\end{align*}
there exists a unique fixed-point $\u^*$ of $T$ in $\B_{\XX_\nu}(\bu,r)$, the closed ball of center $\bu$ and radius $r$ in $\XX_\nu$.
Assume further that $\bw$ (which plays a role in the definition of $A$), is such that
\begin{align*}
\Vert \1 - \bw\ast(1-\bu^{(2)}) \Vert_\nu <1.
\end{align*}
Then $\u^*$ is the unique zero of $F$ in $\B_{\XX_\nu}(\bu,r)$.
\end{theorem}

\subsection{Derivation of the bounds}
\label{sec:bounds_NP}

It once again remains to derive computable estimates $Y$, $Z_1$ and $Z_2$ satisfying the assumptions of Theorem~\ref{th:NK_SKT}. The required computations are still very similar to the one performed in Sections~\ref{sec:boundsPM} and~\ref{sec:bounds_SKT}. 

\subsubsection{The bound $Y$}

$\Vert AF(\bu)\Vert_{\XX_\nu}$ can be computed explicitly using interval arithmetic.

\subsubsection{The bound $Z_1$}

Introducing once more 
\begin{align*}
B=\begin{pmatrix}
B^{(1,1)} & B^{(1,2)} \\
B^{(2,1)} & B^{(2,2)}
\end{pmatrix} = I - ADF(\bu),
\end{align*}
we again consider $Z_1 = \max\left( Z_1^{finite},\ Z_1^{tail} \right)$, where $Z_1^{finite}$ corresponds to
\begin{align*}
 \max_{0\leq n \leq 2N-2} \frac{1}{\xi_n(\nu)} \max\left(\Vert B^{(1,1)}_{(\cdot,n)} \Vert_\nu + \Vert B^{(2,1)}_{(\cdot,n)} \Vert_\nu,\ \Vert B^{(1,2)}_{(\cdot,n)} \Vert_\nu + \Vert B^{(2,2)}_{(\cdot,n)} \Vert_\nu\right),
\end{align*} and
\begin{align}
\label{eq:Z1tail_NP}
Z_1^{tail} &=  \left\vert \left\Vert \begin{pmatrix} \1-\bw\ast(\1-\bu^{(2)}) & \bs - \bw\ast (\bm{\gamma}+\bu^{(1)}) \\ \0 & \0 
\end{pmatrix} \right\Vert_\nu + \frac{1}{(N\pi)^2} \begin{pmatrix} \Vert\bs\Vert_\nu \Vert R'(\bu^{(1)})\Vert_\nu & \0 \\ \Vert R'(\bu^{(1)})\Vert_\nu & \0 
\end{pmatrix}\right\vert_1 \\
& = \max\left[\Vert \1 - \bw\ast(\1 - \bv)\Vert_\nu + \frac{1}{(N\pi)^2} \left( 1 + \Vert \bs\Vert_\nu\right) \Vert R'(\bu))\Vert_\nu, \Vert \bs - \bw\ast (\bm{\gamma}+\bu)) \Vert_\nu\right], \nonumber
\end{align}
bounds
\begin{align*}
\sup_{n\geq 2N-1} \frac{1}{\xi_n(\nu)} \max\left(\Vert B^{(1,1)}_{(\cdot,n)} \Vert_\nu + \Vert B^{(2,1)}_{(\cdot,n)} \Vert_\nu,\ \Vert B^{(1,2)}_{(\cdot,n)} \Vert_\nu + \Vert B^{(2,2)}_{(\cdot,n)} \Vert_\nu\right).
\end{align*}

\begin{remark}
\label{rem:opt}
The obtained $Z_1^{tail}$ estimate shows that the proposed approach of introducing a system of differential-algebraic equations to deal with (at least some types of) non-polynomial diffusion terms is not only valid but also efficient, because we keep the $\frac{1}{N^2}$ scaling. Even the constant is close to optimal, see the discussion in Appendix~\ref{sec:opt}.
\end{remark}

\subsubsection{The bound $Z_2$}

Finally, remembering that $R(u) = \alpha u - \beta u^2 + g$, an estimation similar to the ones performed in the two previous cases leads to
\begin{align*}
Z_2 =  \max\left[\Vert A^{(1,1)} \Vert_{\nu} + \Vert A^{(2,1)} \Vert_{\nu},\, 2\beta\left(\Vert A^{(1,2)} \Vert_{\nu} + \Vert A^{(2,2)} \Vert_{\nu} \right) \right].
\end{align*}

\subsection{Example and results}

In this subsection, we apply the technique we just presented to a specific example, and again consider $R(u) = u - u^2 + g$ (i.e. we take $\alpha=\beta=1$), with $g$ as in Section~\ref{sec:res_PM}, but this time with $\Phi(u) = \frac{u}{\gamma + u}$, and two different values of $\gamma$.

\begin{figure}[h!]
\centering
\begin{subfigure}[t]{0.49\textwidth}
\centering
\includegraphics[scale=0.49]{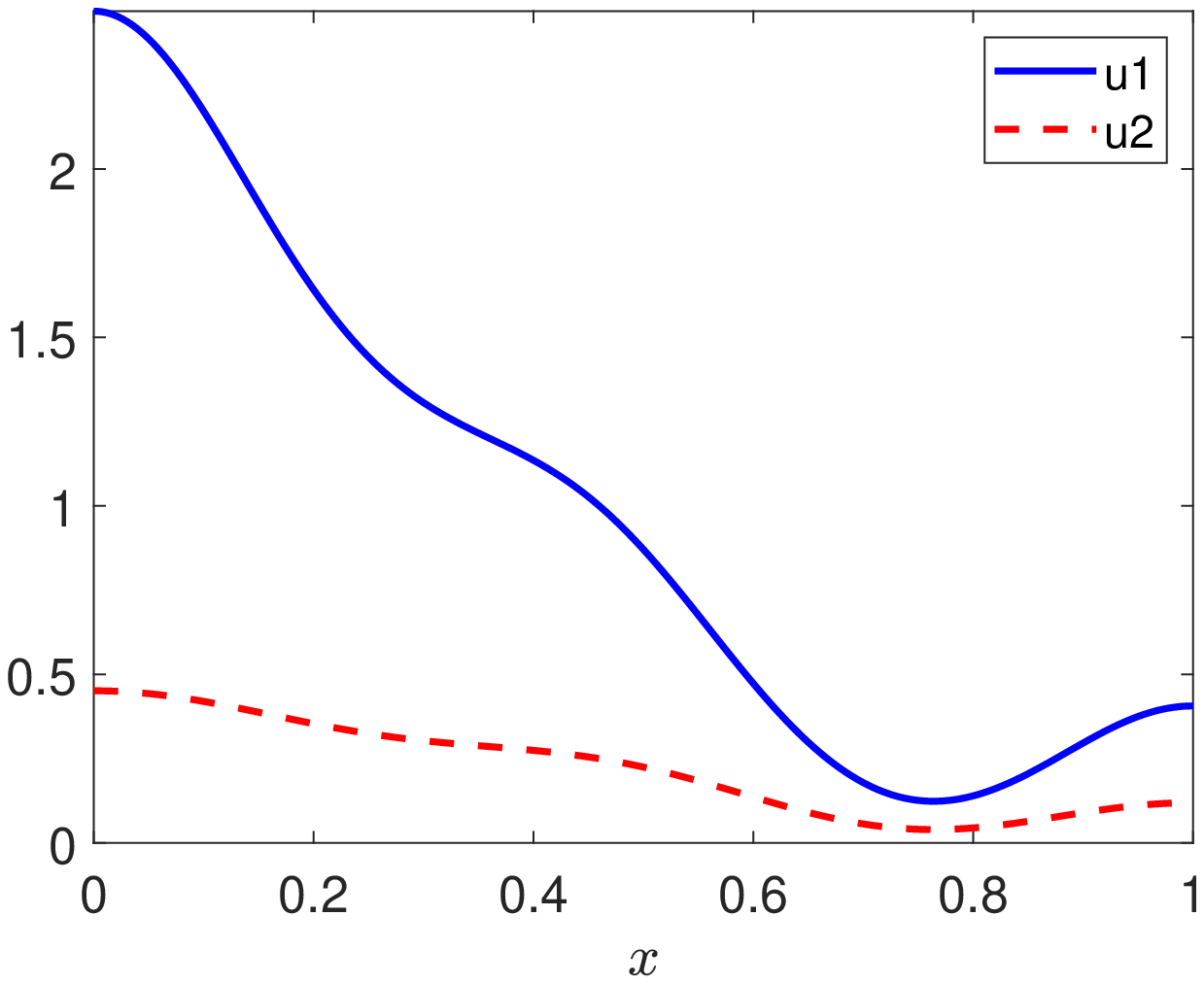}
\caption{$\gamma = 3$.}
\label{fig:NP3}
\end{subfigure}
\begin{subfigure}[t]{0.49\textwidth}
\centering
\includegraphics[scale=0.49]{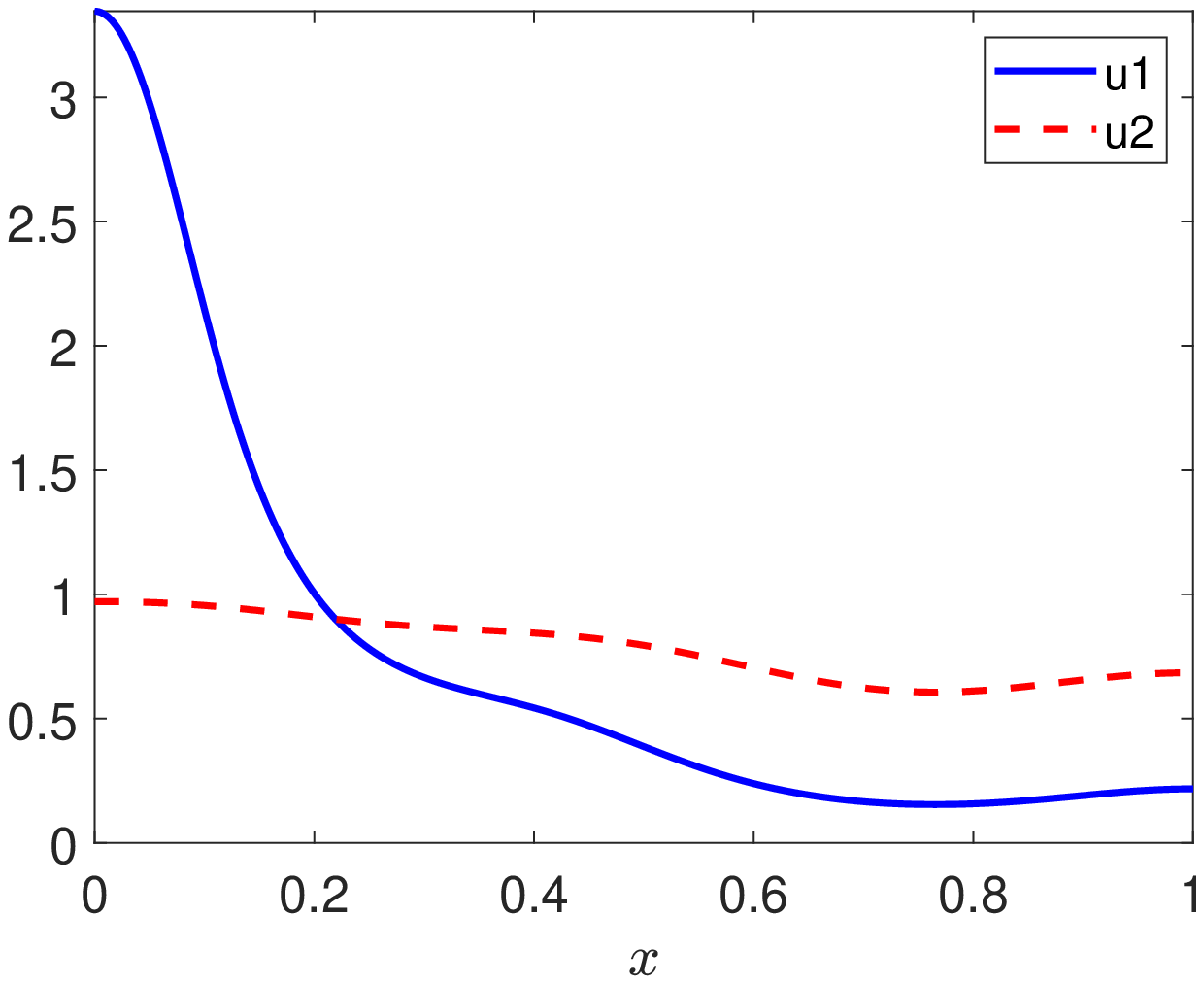}
\caption{$\gamma = 0.1$.}
\label{fig:NP01}
\end{subfigure}
\caption{The two approximate solutions of~\eqref{eq:PMbis} with $\Phi(u) = \frac{u}{\gamma+u}$, $R(u) = u-u^2 +g$ and $g$ as in Figure~\ref{fig:g}, which have been validated in Theorem~\ref{th:NP}. $u_1$ corresponds to the solution $u$ itself, while $u_2$ corresponds to the additional variable $v=\Phi(u)$ introduced in~\eqref{eq:NP}.}.
\label{fig:NP}
\end{figure}

In each case, we computed an approximate solution $\bu$ (see Figure~\ref{fig:NP}), and then validated it using the procedure described in this whole section, which yield the following theorem.

\begin{theorem}
\label{th:NP}
Let $\gamma=3$, and $\bar u $ be the function whose Fourier coefficients $\bu$ can be downloaded at~\cite{Bre21}, and which is represented in Figure~\ref{fig:NP3}. There exists a strong solution $u$ of~\eqref{eq:PMbis}, with $\Phi$ and $R$ as given just above, such that $\left\Vert u-\bar u\right\Vert_{C^0} \leq 1\times 10^{-13}$.

Let $\gamma=0.1$, and $\bar u $ be the function whose Fourier coefficients $\bu$ can be downloaded at~\cite{Bre21}, and which is represented in Figure~\ref{fig:NP01}. There exists a strong solution $u$ of~\eqref{eq:PMbis}, with $\Phi$ and $R$ as given just above, such that $\left\Vert u-\bar u\right\Vert_{C^0} \leq 8\times 10^{-6}$.
\end{theorem}
\begin{proof}
We take $N=50$ and $\nu=1.1$, compute the estimates $Y$, $Z_1$ and $Z_2$ obtained in Section~\ref{sec:bounds_NP} for each approximate solution, and apply Theorem~\ref{th:NK_NP}. 
The computational parts of the proof, can be reproduced by running the Matlab code \texttt{script\_NP.m} from~\cite{Bre21}, together with Intlab~\cite{Rum99}, which provides more details such as the precise value of each bound for each solution.
\end{proof}
Remark~\ref{rem:space} also applies here. To tie back these results to the discussion at the start of this section, notice from Figure~\ref{fig:NP} that, for $\gamma = 3$ we have $\Vert u\Vert_{C^0} <\gamma$, so one could hope to prove that $\Phi(u)= \sum_{k\geq 0} (-1)^k \frac{\u^{k+1}}{\gamma^{k+1}}$ is converging power series, but that is definitely not true for $\gamma=0.1$. In that case, $u$ stays positive so we are still safely away from the pole of $\Phi$, but one would have to use a power series expansion around some non zero value depending on the solution (e.g. the middle of the interval of values taken by $u$).

\appendix
\section*{Appendix}

\section{About the computation of $\bar A$}
\label{sec:barA}

We give here more details about Remark~\ref{rem:barA}, and in particular explain why, no matter how large $N$ is, taking $\bar A = \left(\Pi_N DF(\bu)\Pi_N\right)^{-1}$ is not necessarily suitable. Indeed, let us consider the case described in Section~\ref{sec:PM}, and have a closer look at the coefficient in position $(N,N)$ of $ADF(\bu)$. The coefficients in $A$ and $DF(\bu)$ that are relevant for this computation are highlighted below:
\begin{align*}
\begin{tiny}
\renewcommand\arraystretch{3} 
\newcommand{\ph}{\phantom{\frac{1}{N^2} \bw_1}}
\left(\begin{array}{cccccc}
 & & & \multicolumn{1}{|c}{ }  &   &   \\ 
 & \bar A & & \multicolumn{1}{|c}{ } &  &  \\ 
 \ph & \ph & \ph & \multicolumn{1}{|c}{\frac{-\bw_1}{(N\pi)^2} } & \frac{-\bw_2}{((N+1)\pi)^2}   & \ldots \\ \cline{1-3}
  &  &  &  &  &   \\
  &  &  &  &  &   \\
  &  &  &  &  &   \\
\end{array}\right)
\left(\begin{array}{ccccc}
 & & & \multicolumn{1}{|c}{ }  &    \\ 
 & \Pi_N DF(\bu)\Pi_N & & \multicolumn{1}{|c}{ } &  \\ 
 \ph & \ph & \ph & \multicolumn{1}{|c}{ \ph } &  \ph  \\ \cline{1-3}
  &  & -N^2\Phi'_1(\bu) &  &    \\
  &  & -(N+1)^2\Phi'_2(\bu) &  &    \\
  &  & \vdots &  &    \\
\end{array}\right).
\end{tiny}
\end{align*}
There should also be terms corresponding to $R'(\bu)$ in $DF(\bu)$, but they do not have a factor $N^2$ in front, hence for $N$ large enough we can neglect them. Therefore we have
\begin{align*}
\left( ADF(\bu)\right)_{N,N} \approx \left( \bar A \ \Pi_N DF(\bu)\Pi_N\right)_{N,N} + \sum_{n\geq 1} \bw_n \Phi'_n(\bu).
\end{align*}
Now, if we were to take $A\approx \left(\Pi_N DF(\bu)\Pi_N\right)^{-1}$, then we would have in particular 
\begin{align*}
\left( \bar A \ \Pi_N DF(\bu)\Pi_N\right)_{N,N} \approx 1
\end{align*}
and therefore
\begin{align*}
\left( ADF(\bu)\right)_{N,N} \approx 1 + \sum_{n\geq 1} \bw_n \Phi'_n(\bu).
\end{align*}
However, $\sum_{n\geq 1} \bw_n \Phi'_n(\bu)$ has no reason to be small in general, so $\left( ADF(\bu)\right)_{N,N}$ could be significantly different from $1$, which would prevent us from satisfying~\eqref{eq:Z1l1}.

\section{About the $Z_1^{tail}$ bound in the non-polynomial case}
\label{sec:opt}

We make some comments about the quality of the bound~\eqref{eq:Z1tail_NP}, which in some sense measures how well adapted the transformation of~\eqref{eq:PMbis} into~\eqref{eq:NP} and the choice of $A$ in~\eqref{eq:A_NP} are for computer-assisted proofs, and explain why we claim in Remark~\ref{rem:opt} the output is close to optimal. 

In order to get a comparison, notice that for this example $\Phi$ can be inverted by hand, and therefore~\eqref{eq:PMbis} can easily be rewritten as an equation with linear diffusion
\begin{align*}
\left\{
\begin{aligned}
&\Delta v + \tilde R(v) = 0 \qquad & \text{on }(0,1), \\
&\frac{\partial v}{\partial n } = 0 \qquad & \text{on } \{0,1\},
\end{aligned}
\right.
\end{align*}
where $v=\Phi(u)$ and
\begin{align*}
\tilde R(v) = R\left(\Phi^{-1} (v)\right) = R\left( \frac{\gamma v}{1-v}\right).
\end{align*}
In this case which is standard since there is no nonlinear diffusion, the $Z_1^{tail}$ estimate would be equal to
\begin{align*}
\frac{1}{(N\pi)^2} \Vert \tilde R'(\bv) \Vert_\nu & = \frac{1}{(N\pi)^2} \left\Vert \frac{\gamma}{(1-\bv)^2} R'(\bu) \right\Vert_\nu.
\end{align*}
However, since $\bw\approx \frac{1}{1-\bv}$, $\bm\gamma+\bu\approx\bm\gamma\bw$ and $\bs\approx \bw\ast(\bm{\gamma}+\bu) \approx \gamma\bw^2$, hence the above $Z_1^{tail}$ estimate would be roughly equal to
\begin{align*}
\frac{1}{(N\pi)^2} \left\Vert \bs R'(\bu) \right\Vert_\nu.
\end{align*}
The estimate obtained in~\eqref{eq:Z1tail_NP} is worst, but not dramatically so. Indeed, in practice the first part of~\eqref{eq:Z1tail_NP} can always be made small enough not to matter, and essentially~\eqref{eq:Z1tail_NP} reduces to
\begin{align*}
\frac{1}{(N\pi)^2} (1+\left\Vert \bs\right\Vert_\nu) \left\Vert R'(\bu) \right\Vert_\nu.
\end{align*} 
The extra $1$ could even be removed by choosing an appropriately weighted norm on the product space $\ell^1_\nu\times\ell^1_\nu$ (see e.g. the discussion about the weights in~\cite{BreKue19}). Therefore the only thing that we lost, at least concerning the crucial $Z_1^{tail}$ estimate, by going through a system of differential-algebraic equation is that $\left\Vert \bs R'(\bu) \right\Vert_\nu$ is replaced by $\left\Vert \bs\right\Vert_\nu \left\Vert R'(\bu) \right\Vert_\nu$, which in most cases is perfectly acceptable. This shows that, in more complex situations where $\Phi$ cannot be explicitly inverted, the approached proposed in this paper provides a good alternative.

\bibliographystyle{abbrv}
\bibliography{bibfile}

\end{document}